\input ./style/arxiv-general.cfg
\documentclass[aop,MSNbibl,seceqn,dvips]{arximspdf}
\makeatletter
   \@ifpackageloaded{graphicx}{}{\usepackage{graphicx}}
\makeatother
\usepackage{mathrsfs}

\doi{10.1214/14-AOP957}
\volume{43}
\issue{6}
\pubyear{2015}
\firstpage{3177}
\lastpage{3215}
\docsubty{FLA}

\makeatletter

\def\afrac#1#2{#1/(#2)}

\newcommand{\rrvert}{\vert}
\newcommand{\rrVert}{\Vert}
\newcommand{\llvert}{\vert}
\newcommand{\llVert}{\Vert}
\renewcommand{\mid}{|}
\newtheorem{theorem}{Theorem}[section]
\newtheorem{proposition}[theorem]{Proposition}
\newtheorem{lemma}[theorem]{Lemma}
\newproclaim{remark}[theorem]{Remark}
\newproclaim{definition}[theorem]{Definition}
\newproclaim{claim}{Claim}
\newcommand{\R}{\mathbb{R}}
\newcommand{\D}{\mathcal{D}}
\newcommand{\bQ}{\mathbf{Q}}
\newcommand{\Z}{\mathbb{Z}}
\newcommand{\PATH}{\zeta}
\def\eps{\varepsilon}
\def\P{\mathbf{P}}
\def\E{{\mathbb E}}
\def\dc{d_{\mathrm{sup}}}
\def
\prob{\mathbb{P}}
\def\TUBE{\mathcal{T}}
\def\HH{\mathscr{H}} 
\def\K{\mathcal{K}}
\def\calW{\mathcal{W}}
\def\Cr{\mathsf{Cr}}
\def\PATH{\zeta}
\def\R{\mathbb{R}}
\def\cW{\mathcal{W}}
\def\cV{\mathcal{V}}
\def\cS{\mathcal{S}}
\def\HH{\mathscr{H}}
\def\cE{\mathcal{E}}
\def\cD{\mathcal{D}}
\def\cB{\mathcal{B}}
\def\boxT{\mathcal{T}_{\square}}
\def\triT{\mathcal{T}_{\triangle}}
\def\tr{\triangle}
\def\LL{\mathbb{L}}
\newcommand{\osc}{\operatorname{osc}}
\newcommand{\margin}[1]{} 
\makeatother

\begin{document}
\begin{frontmatter}

\title{Coalescing Brownian flows: A new approach}
\runtitle{Coalescing Brownian flows: A new approach}

\begin{aug}
\author[A]{\fnms{Nathana\"el}~\snm{Berestycki}\corref{}\thanksref{T1}\ead[label=e1]{nberestycki@gmail.com}},
\author[B]{\fnms{Christophe}~\snm{Garban}\thanksref{T2}\ead[label=e2]{garban@math.univ-lyon1.fr}}
\and
\author[C]{\fnms{Arnab}~\snm{Sen}\ead[label=e3]{arnab@umn.edu}}
\runauthor{N. Berestycki, C. Garban and A. Sen}
\affiliation{University of Cambridge, Universit\'e Lyon 1 and University of
Minnesota}
\address[A]{N. Berestycki\\
Statistical Laboratory\\
DPMMS\\
University of Cambridge\\
Wilberforce Rd.\\
Cambridge, CB3 0WB\\
United Kingdom\\
\printead{e1}}
\address[B]{C. Garban\\
Institut Camille Jordan\\
Universit\'e Lyon 1\\
43 bd du 11 novembre 1918\\
69622 Villeurbanne cedex\\
France\\
\printead{e2}}
\address[C]{A. Sen\\
Department of Mathematics\\
University of Minnesota\\
127 Vincent Hall\\
206 Church St. SE\\
Minneapolis, Minnesota 55455\\
USA\\
\printead{e3}}
\end{aug}
\thankstext{T1}{Supported in part by EPSRC Grants EP/GO55068/1 and
EP/I03372X/1.}
\thankstext{T2}{Supported in part by the ANR Grant MAC2 10-BLAN-0123.}

%
\received{\smonth{11} \syear{2013}}
%
\revised{\smonth{7} \syear{2014}}

%
\begin{abstract}
The coalescing Brownian flow on $\mathbb{R}$ is a process which was
introduced by
Arratia [Coalescing Brownian motions on the line (1979) Univ. Wisconsin, Madison]
and T\'oth and Werner [\textit{Probab. Theory Related Fields} \textbf{111} (1998) 375--452],
and which formally corresponds to
starting coalescing
Brownian motions from every space--time point. We provide a new state
space and topology
for this process and obtain an invariance principle for coalescing
random walks. This
result holds under a finite variance assumption and is thus optimal. In
previous works by
Fontes et al. [\textit{Ann. Probab.} \textbf{32} (2004) 2857--2883],
Newman et al. [\textit{Electron. J. Probab.} \textbf{10} (2005) 21--60],
the topology and
state-space required a moment
of order $3-\varepsilon$ for this convergence to hold. The proof relies
crucially on recent work
of Schramm and Smirnov on scaling limits of critical percolation in the
plane. Our
approach is sufficiently simple that we can handle substantially more
complicated
coalescing flows with little extra work---in particular similar
results are obtained
in the case of coalescing Brownian motions on the Sierpinski gasket.
This is the first
such result where the limiting paths do not enjoy the \emph{noncrossing property}.
\end{abstract}

%
\end{frontmatter}

\section{Introduction}

\subsection{Motivation}
The coalescing Brownian flow or Arratia's flow was first introduced and
studied by Arratia \cite{arratia79,arratia81}, as a limiting object
describing the large-scale behaviour of the one-dimensional voter
model. Informally, this process consists of particles that perform
independent coalescing Brownian motions, starting from every space time
$(x, t)$. By \emph{independent coalescing Brownian motions}, we mean
two paths which are independent Brownian motions until the first time
they meet, and which subsequently continue as one single Brownian motion.

Since then, the Brownian web has been conjectured or proved to describe
scaling limits in a large number of seemingly disconnected models: let
us mention in particular the works of T\'oth and Werner \cite{toth98}
in connection with true self-repelling motion, Coletti et al. \cite
{coletti09} in connection with a drainage network model,
Sarkar and Sun \cite{sarkar13} in connection with oriented
percolation, Norris and Turner \cite{norris12} in connection with
Hastings--Levitov planar aggregation models.

It is not a priori easy to turn the informal description in the first
paragraph into a rigorous mathematical object. The difficulty lies in
the fact that there are uncountably many starting points. Curiously,
foundations for constructing the Brownian web as a random variable in a
``nice'' space and studying convergence of discrete objects to the
Brownian web were only laid down recently, in a series of papers by
Fontes et al. \cite{44444}. 


\subsection{Main results}

Our object in this paper is threefold: \margin{C: used to be twofold,
so I changed a bit the order here and insisted from the beginning on
the Sierpinski flow. I changed a bit the wording here and there, please check.}
\begin{longlist}[2.]
\item[1.] First, we provide an alternative state-space and topology for
coalescing flows such as Arratia's flow (Theorem~\ref{tBW}). Inspired
by the setup introduced by Schramm and Smirnov in \cite
{SchrammSmirnov11} for critical percolation, we call this space the
\textit{Schramm--Smirnov space of coalescing flows}
$\HH$. See Definition~\ref{dspace}. As the reader will see, our
setup has the great advantage that it makes the proof of convergence of
discrete objects to the limiting Brownian flow surprisingly simple.

\item[2.] We then prove an invariance principle (Theorem~\ref
{tscalinglimit}) for the convergence of scaled coalescing random walks
on $\mathbb{Z}$ toward Arratia's flow under an optimal finite variance
assumption on the random walk.
Note that in previous works, the topological setup was different and
establishing tightness already required a nontrivial proof. In
particular, in \cite{bmsv06}, it is shown that in order to obtain a
tightness criterion, a $3+\eps$ finite moment is sufficient, while
existence of a $3-\eps$ moment is needed. We prove in our present
setting that a finite variance is both necessary and sufficient.
\margin{C: I changed according to latest's Arnab comment. Please check
Nathanael if you are
ok with this.}


\item[3.]
Finally, we illustrate the simplicity and flexibility of our approach
by showing similar results for coalescing flows where the underlying
geometry is substantially more complicated.
We focus in particular on coalescing random walks on the discrete
infinite Sierpinski gasket, and prove that this can be rescaled to a
coalescing Brownian flow on the continuous Sierpinski gasket (Theorem
\ref{tscalinglimitg}). On the real line any two continuous paths
cannot cross without hitting each other.
This is an obvious topological fact in dimension one which underlies
Arratia's original approach and much of the work on the subject. But
this property is absent for the Sierpinski Gasket since one of the
paths can go into some other ``triangle'' and come back later at a
suitable time.
\end{longlist}



At the heart of our approach is the groundbreaking work of
Schramm and Smirnov \cite{SchrammSmirnov11} on scaling limits of
critical percolation in the plane.
We borrow directly from their work by adapting the quad-crossing or
Schramm--Smirnov state-space and topology to the study of coalescing
flows. As in the percolation case, the state space is compact, so only
uniqueness of subsequential limits has to be established. Thus, one of
the main advantages of our approach is that very few estimates are
needed in order to establish this convergence. The bulk of the work is
in some sense to construct the limiting object in the Schramm--Smirnov
space. Convergence of the discrete object to the continuous one then
follows rather simply by an argument based on ``uniform coming down
from infinity'' (Proposition~\ref{propcomingdownfrominf}), which
is the main technical ingredient.

 In the companion paper \cite{SBN}, \margin{C: added this
reference} we will illustrate further the extent and breadth of the
analogy between coalescing flows and critical percolation, by showing
that the Brownian webs in this paper (including, say, on the Sierpinski
gasket) are all examples of \emph{black noises}.

\subsection{Relation to previous work}

As mentioned earlier, one of the main results in this paper (Theorem
\ref{tscalinglimit}) is an invariance principle for coalescing random
walk under a second moment assumption. This contrasts sharply with the
approach developed by Fontes et al. \cite{44444}: indeed, using that
topology, a series of works culminating with Newman et al. \cite
{nrs05} showed that such an invariance principle holds a under a fifth
moment assumption. Subsequently, Belhaouari et al. \cite
{bmsv06} lowered this to a $3+\eps$ moment assumption and,
surprisingly, showed that this was in fact also necessary, in the sense
that the convergence does not hold if the random walk's $3-\eps$
moments are not finite for some $\eps>0$.

To understand better the difference between the two approaches, it is
useful to\margin{added to, please check} dwell further on the analogy
with percolation.
In that context, it is also a nontrivial task to build a ``nice''
state-space and topology for taking scaling limits (as the mesh size
tends to 0). One can study the collection of all individual contour
interfaces of clusters (an approach originating in the work of
Aizenmann and Buchard \cite{aizenman99} and culminating in the work of
Camia and Newman \cite{camia06}), or one can ask about macroscopic
connectivity properties. This is the viewpoint taken in the work of
Schramm and Smirnov \cite{SchrammSmirnov11} which motivates our
approach. This latter topological framework introduced by Schramm and
Smirnov was recently used and extended in \cite{GPS2a} and \cite
{GPS2b} in order to prove that \textit{near-critical percolation} on the
triangular grid has a (massive) scaling limit. \margin{C: added this
advertisement:) }

Some alternative approaches to the construction of Arratia's
flow have been proposed which we briefly review. One of the problem for
formulating a state-space for Arratia's flow is that it is not in the
strict sense a flow, for example, it may well be the case that from a point
$x\in\R$ there is more than one trajectory coming out of it. To go
around this problem, Norris and Turner \cite{norris12} considered the
space of \emph{weak flows}, which replaces the notion of a flow $(\phi
_{st})_{s<t}$ by a pair $(\phi^{\pm}_{st})_{s<t}$ satisfying some
compatibility relations. Essentially, $\phi^-$ is the left-continuous
version of the flow and $\phi^+$ the right-continuous version.
Arratia's flow is then identified as the unique random variable on this
space such that its restriction to finitely many points forms a
coalescing Brownian motion started from that set. This approach is very
elegant but has the drawback that it relies crucially on the
noncrossing property.


Le Jan and Raimond, in a series of papers, \cite{lejan04a,lejan04b},
adopted a
different point of view on the question. They viewed Arratia's flow as
a (random) flow of maps, which are essentially consistent systems of
$n$-point Markovian systems, describing the law of the motion of $n$
indivisible points. While this approach is in principle very general
(and in particular, does not rely on the noncrossing property), it is
not well-suited to the questions of taking a scaling limit of some
discrete flow to its continuous counterpart. This is because the
question of scaling limits would have to be approached through the
finite-dimensional distributions of the flow. See, for instance, \cite{LeJanLemaire}
for an example and see \cite{LeJanICM} for a good survey of this approach.

\subsection{Organisation of the paper}

The paper will be divided as follows:
\begin{itemize}
\item
In Section~\ref{Stopology}, we describe the setup used in this paper,
that is, the Schramm--Smirnov space $\HH$ (see Definition~\ref
{dspace}) and its topology. We also explain how to view a compact set
of paths as an element of the state space $\HH$, and give a convenient
criterion for convergence in terms of tube-crossing probabilities.

\item In Section~\ref{SBW}, we give a construction and
characterisation of the Brownian web as a random variable on $\HH$
(Theorem~\ref{tBW}). This is immediately followed in Section~\ref
{Scharac} by a result characterising further the Brownian web. This is
not needed for the rest of the paper and may be skipped by a reader who
only wishes to read the invariance principle of Section~\ref{Sscalinglimit}.

\item In Section~\ref{Sscalinglimit}, we state and prove the
invariance principle showing that rescaled coalescing random walks
converge to the Brownian web (Theorem~\ref{tscalinglimit}).

\item In Section~\ref{Sgasket}, we extend our setup to study simple
coalescing random walks on the Sierpinski gasket.

\item In the final Section~\ref{Sgasket2}, we state and prove an
invariance principle showing that rescaled coalescing random walks on
the Sierpinski graph converge to the coalescing flow on the Sierpinski
gasket constructed in Section~\ref{Sgasket}. \margin{C: added this
item which was missing.}

\item In the \hyperref[app]{Appendix}, we make the link with the state-space and
topology of Fontes et al. In particular, we show that our topology is coarser.
\end{itemize}



\section{Schramm--Smirnov topology and flows}\label{Stopology}

\subsection{The space of tubes}

We start by introducing the notion of tubes (which replaces the notion
of quads in the case of planar percolation). Fix $d\ge1$;
we will consider in this section coalescing flows for which the
individual paths are naturally embedded into $\R^d$. We will denote a
generic space--time point $z=(x,t)$, where $x\in\R^d$ represents the
spatial coordinate and $t$ is the time-coordinate.

%
\begin{figure}

\includegraphics{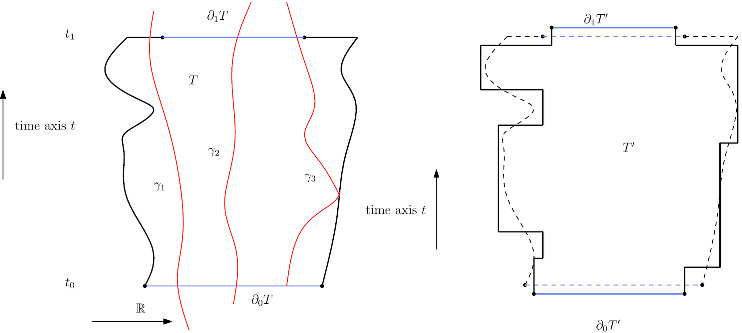}

\caption{On the left, a tube $T$ with lower and upper face $\partial
_0 T$ and $\partial_1 T$. The tube $T$ is traversed by paths $\gamma
_2$ and $\gamma_3$ but not by $\gamma_1$. On the right, a tube $T'$
from the class
$ \TUBE_\square$ of Definition \protect\ref{Dnicetubes}
approximating a tube $T$, shown in dotted line.}\label{ftube}
\end{figure}

\begin{definition}[(Tube)]\label{dtube}
A tube $T$ is a triplet $([T], \partial_0 T, \partial_1 T):=\break 
(\varphi([0, 1]^{d+1}), \varphi([0,1]^d\times\{0\}), \varphi
([0,1]^d\times\{1\}) )$ where $\varphi\dvtx  [0,1]^d \times[0,1] \to\break 
\varphi( [0,1]^d \times[0,1]) \subset\R^d \times\R$ is a homeomorphism
such that $\varphi([0,1]^d\times\{0\})$
and $\varphi([0,1]^d \times\{1\})$ are subsets of $\R^d \times\{
t_0\}$ and $\R^d \times\{t_1\}$, respectively, for some $t_0<t_1$.
Furthermore, we require that $[T]=\varphi([0,1]^{d+1})$ is included in
$\R^d \times[t_0,t_1]$. \margin{Added this easier constraint
suggested by Arnab to fix the issue raised before.}
Informally, $T$ is a topological cube along with a distinct pair of
opposite faces which are both orthogonal in $\R^{d+1}$ to the time axis.
We call $t_0$ the {\em start time} of $T$ and $t_1$ the {\em end time}
of $T$. The sets $\partial_0 T$ and $\partial_1 T$ are called the
{\em lower face} and the {\em upper face} of $T$, respectively.
See Figure~\ref{ftube} in the case $d=1$.
\end{definition}

%
%

\begin{definition}[(Metric space of tubes)]\label{dTUBE}
The space of all such tubes, denoted by $\TUBE$, can be equipped with
the following metric:
\[
d_\TUBE(T_1,T_2):= d_{\mathrm{Haus}}
\bigl([T_1], [T_2]\bigr) + d_{\mathrm
{Haus}} (
\partial_0 T_1, \partial_0 T_2) +
d_{\mathrm
{Haus}}(\partial_1 T_1, \partial_1
T_2),
\]
where $d_{\mathrm{Haus}}$ is the usual Hausdorff metric on the compact
subsets of $\R^{d+1}$, which is given by
$d_{\mathrm{Haus}}(A, B) = \max( \sup_{ b \in B} \inf_{a \in A} \llVert
a-b\rrVert  _2, \sup_{ a \in A} \inf_{ b \in B} \llVert   a-b\rrVert  _2 )$. It is easy to
see that $(\TUBE, d_\TUBE)$ is separable (see, e.g., Section~\ref{dktubes}).
\end{definition}

\subsection{The space of coalescing flows}

First, we introduce the notion of \emph{crossing} or \emph
{traversing} for a tube by a continuous path in $\R^d$ with a
specified starting time. A continuous path $\gamma$ with starting time
$t$ is just a continuous map $\gamma\dvtx [t, \infty) \to\R^d$ and is
denoted by $(\gamma, t)$. We will say that a tube $T$ is \emph
{crossed} or \emph{traversed} by $(\gamma, t)$ if $t \le t_0$,
$(\gamma(t_0), t_0) \in\partial_0 T$, $(\gamma(t_1), t_1) \in
\partial_1 T$ and $(\gamma(s), s) \in[T]$ for all $s \in(t_0,
t_1)$, where $t_0$ and $t_1$ are the start time and the end time of
$T$, respectively. Informally, this means that the trajectory of the
path $\gamma$ enters the tube $T$ through the lower face $\partial_0
T$ and stays in it throughout until it exits $T$ through its upper face
$\partial_1 T$. Given any (finite or infinite) collection of such
continuous paths with starting times, we can associate the subset of
$\TUBE$ which consists of all the tubes which are crossed by at least
one path from the collection. For example, in the case $d=1$, for a
family of continuous (random) paths which are trajectories of some
coalescing Brownian particles starting from different space--time points
in $\R^2$, we can naturally obtain via the above association an
element of $\{0,1\}^{\TUBE}$. This contains a lot of information about
the coalescing paths.
This suggests taking $\{0,1\}^{\TUBE}$ as a state-space for the
Brownian web.
But as indicated in \cite{SchrammSmirnov11}, this turns out to be too
big and unwieldy for our purpose. Instead, note that a subset of $\TUBE
$ consisting of all the tubes which are crossed by some collection of
paths cannot be arbitrary: some compatibility conditions must be
satisfied. The following partial order on tubes, which mirrors that of
\cite{SchrammSmirnov11} will help to clarify the structure.

\begin{figure}

\includegraphics{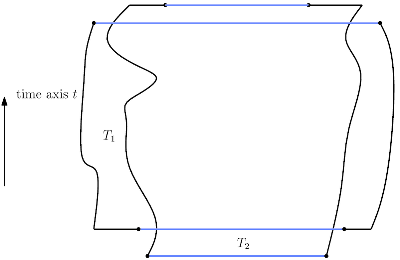}

\caption{Two tubes such that $T_1 < T_2$.}\label{fporder}
\end{figure}

\begin{definition}[(Partial orders on tubes)]\label{dpartialorder}
If $T_1, T_2 \in\TUBE$ are two tubes, we will say that:
\begin{longlist}[(a)]
\item[(a)] $T_1 \le T_2$ if whenever $T_2$ is traversed by any path
$\gamma$, $T_1$ is also traversed by~$\gamma$. See Figure~\ref{fporder}.
\item[(b)] $T_1 < T_2$ if there are open neighborhoods $U_i$ of $T_i$
in $(\TUBE, d_\TUBE)$ such that $T_1' \le T'_2$ holds
for any $T'_i \in U_i$, $i=1,2$.
\end{longlist}
\end{definition}

%

The definition below (following again \cite{SchrammSmirnov11}) singles
out a particular class of subsets of $\TUBE$
by adding the compatibility constraint on the tubes which are
considered to be traversed.

\begin{definition}[(Hereditary subsets of $\TUBE$)]\label{dhereditary}
A subset $ S\subset\TUBE$ is called {\em hereditary} if whenever
$T\in S$ and $T' \in\TUBE$ such that $T' < T$,
we also have $T' \in S$.
\end{definition}

We are now ready to define the crucial definition for our state-space:

\begin{definition}[(Space of coalescing flows)]\label{dspace}
Let $\HH$ be the space of all {\em closed} hereditary subsets of
$\TUBE$. We call $\HH$ the Schramm--Smirnov space of coalescing flows
or simply the space of coalescing flows.
\margin{C: I changed the name, I did not like the term percolative
flows. Please check if you are fine.}
\end{definition}

\begin{remark}
The word \emph{closed} in the above definition refers to the Hausdorff
metric on tubes of Definition~\ref{dTUBE}. Essentially, it is a way
of setting the following convention: if $T_n \to T$ in $d_\TUBE$ and
each $T_n$ is crossed by a family of continuous paths, then we declare
$T$ to be crossed by that family.
Another way to think about it is that it shows how to decide whether a
tube is crossed in the ambiguous case when a path touches the boundary
of the tube, like the path $\gamma_3$ in Figure~\ref{ftube}. Such a
tube is declared crossed.
\end{remark}


\subsection{The tube-topology}

Having defined our state-space $\HH$ we need a convenient topology for it.
Following \cite{SchrammSmirnov11}, we define a topology on $\HH$
generated by the following subbase (their finite intersections form a
basis for the topology)
%
\begin{eqnarray}
&\displaystyle \{ H \in\HH\dvtx  H \cap U \ne\varnothing\}, \qquad U \subseteq\TUBE \mbox{ topologically open},& \label{dclosed1}
\\
&\displaystyle \{ H \in\HH\dvtx  T \notin H \}, \qquad T \in\TUBE.&\label{dclosed2}
\end{eqnarray}
We denote $\{ H \in\HH\dvtx  T \in H \}$ and $ \{ H \in\HH\dvtx  H \cap U \ne
\varnothing\}$ by $\boxminus_T$ and $\boxminus^U$, respectively.

Intuitively speaking, this choice of open sets for the tube-topology
leads to the following desirable properties: if a tube $T$ is crossed
by a collection of paths, and if those paths undergo a slight
`perturbation', then some tube in the neighborhood of $T$ will still be
crossed [corresponding to (\ref{dclosed1})]. Also, if $T$ is not
crossed by a collection of paths, then a small perturbation in those
paths cannot result in a crossing of $T$ [corresponding to (\ref{dclosed2})].

The following fundamental result is due to Schramm and Smirnov.

\begin{theorem}\label{tcompact}
The Schramm--Smirnov space $\HH$ equipped with the tube-topology has
the following properties:
\begin{longlist}[2.]
\item[1.] $\HH$ is separable, metrizable, and hence a Polish space.
\item[2.] $\HH$ is compact.
\item[3.]  For any dense set $\mathcal{S}$ of $\TUBE$, the $\sigma$-field
generated by $\{ \boxminus_T\dvtx  T \in\mathcal{S}\}$ is the Borel
$\sigma$-field on $\HH$.
\end{longlist}
\end{theorem}

The above theorem is an immediate consequence of a more general result
(precisely Theorem~3.10) of \cite{SchrammSmirnov11} which works for
the space of closed hereditary subsets of any second-countable
topological space $X$ with a partial ordering $\prec$ such that the
ordering satisfies
%
\begin{equation}
\label{ttopo1} \bigl\{(x, y) \in X^2\dvtx  x \prec y \bigr\}\qquad
\mbox{is a topologically open subset of } X^2,
\end{equation}
and
%
\begin{equation}
\label{ttopo2} \forall x \in X, \qquad x \in\overline{\{ y \in X\dvtx  y \prec x\}}.
\end{equation}
For us, $X = \TUBE$ with the partial order $<$ in which case
conditions (\ref{ttopo1}) and (\ref{ttopo2}) are easy to verify.

The reader may find close similarities between the tube-topology and
the standard Fell topology (or the topology for closed convergence)
defined on $\operatorname{Cld}(\TUBE)$, the space of closed subsets of
$\TUBE$. The Fell topology has many desirable properties including
compactness [which is lacked by the topology induced by the Hausdorff
metric on $\operatorname{Cld}(\TUBE)$]. The fact is using (\ref{ttopo1})
and (\ref{ttopo2}), one can prove that $\HH$ is closed subset of
$\operatorname{Cld}(\TUBE)$ under the Fell topology, and moreover, the
induced topology on $\HH$ is the same as the tube-topology.

Later we will need the next lemma which says that a monotone (with
respect to set inclusion) sequence in $\HH$ always converges.

\begin{lemma}\label{lmonotone}
Let $H_1 \subseteq H_2\subseteq\cdots$ be a nondecreasing sequence
in $\HH$. Then there exists $H \in\HH$ such that $d_{\HH}(H_n, H)
\to0$. Furthermore, $H= \overline{\bigcup_{n=1}^\infty H_n }$, the
closure being taken with respect to $d_{\TUBE}$.
\end{lemma}


\begin{pf} 
Since $\HH$ is compact, $H_k$ has subsequential limits. Let $H, H'$ be
two such limits, taken, respectively, along the subsequences $(n_k)$ and
$(m_j)$. Let us show first that $H_n \subseteq H$ for all $n \ge1$. If
$T$ is a tube such that $T\in H_n$ for some $n$ then by monotonicity $T
\in H_{n_k}$ for all $k$ sufficiently large. In other words, $H_{n_k}
\in\boxminus_T$ for all sufficiently large $k$. Since $ \boxminus_T$
is a closed set in $(\HH, d_\HH)$, we have $H = \lim_{k} H_{n_k}
\in\boxminus_T$. Hence, $T \in H$ as well.


Now let us show that $H' \subseteq H$. Fix $T \in H'$. By definition of
the basic open sets given in (\ref{dclosed1}), for all $\delta>0$
there exists $T_\delta\in\TUBE$ such that $d_\TUBE(T, T_\delta) <
\delta$ and $T_\delta\in H_{m_j}$ for all sufficiently large $j$.
Since $H_{m_j} \subseteq H$ by the above observation, $T_\delta\in H$
for all $\delta$. But now since $H$ is a closed subset of $\TUBE$ and
$d_\TUBE( T_\delta, T) \to0$ as $\delta\to0$, it follows that
$T\in H$ as well. Hence, $H' \subseteq H$. By symmetry $H = H'$ and
this proves the uniqueness of the subsequential limit.

Now since $H_k \subseteq H$ for all $k \ge1$ and $H$ is closed, we
have $\overline{\bigcup_{k=1}^\infty H_k } \subseteq H$. We have also
shown that given $T \in H$ there exists a sequence of tubes $T_\delta
\in\bigcup_{k=1}^\infty H_k $ such that $d_\TUBE( T_\delta, T) \to
0$ as $\delta\to0$. This shows that $H \subseteq\overline{\bigcup_{k=1}^\infty H_k }$ and the proof of the lemma is complete.
\end{pf}

\subsection{Coalescing paths as an element of the space $\mathscr{H}$}

Let $C[t_0]$ denote the set of continuous functions from $[t_0, \infty
)$ to $\R^d$. Define
\[
\Pi= \bigcup_{ t_0 \in\R} C[t_0] \times\{
t_0\},
\]
where $(\gamma, t_0) \in\Pi$ represents a path $\gamma$ in $\R^d$
which starts from time $t_0$. For $(\gamma, t_0) \in\Pi$, let $\hat
\gamma$ be the continuous function that extends $\gamma$ to entire
time $\R$ by setting $\hat\gamma(t) = \gamma(t_0)$ for all $t< t_0$.
We then define a distance $\varrho$ on $\Pi$ by
%
\begin{equation}
\label{pathmetric} \varrho\bigl((\gamma_1, t_1), (
\gamma_2, \gamma_2)\bigr) = \dc(\hat\gamma _1,
\hat\gamma_2) + \llvert t_1 - t_2\rrvert,
\end{equation}
where $\dc(\hat\gamma_1,\hat\gamma_2) = \sum_{ k=1}^\infty2^{-k}
\sup_{ t \in[-k, k]} \min(\llVert  \hat\gamma_1(t) - \hat\gamma_2(t)\rrVert
_2, 1)$. It is easy to check that $(\Pi, \varrho)$ is a complete
separable metric space.

If $\PATH\subseteq\Pi$ is any collection of continuous functions
(with starting times), one can naturally associate to it a subset
$\Cr(\PATH)$ of $\TUBE$ consisting of all tubes traversed by~$\PATH
$. It is clear that $\Cr(\PATH)$ is hereditary for any $\PATH
\subseteq\Pi$. To associate to $\PATH$ a genuine member of $\HH$ we
need to check that $\Cr(\PATH)$ is closed in $\TUBE$ so that $\Cr
(\PATH) \in\HH$. The following straightforward lemma says that this
happens at least if $\PATH$ is compact.

%
%

\begin{lemma} \label{compactimpliescont}
If $\PATH\subseteq\Pi$ is compact, then $\Cr(\PATH) \in\HH$.
\end{lemma}

\begin{pf}
We need to show that if $T_n \in\Cr(\PATH)$ and $d_\TUBE(T_n, T)
\to0$, then $T \in\Cr(\PATH)$.
Suppose that $T_n$ is traversed by $(\gamma_n, \tau_n) \in\PATH$.
Since $\PATH$ is compact, we can find a subsequence $\{n_k\}$ such
that $(\gamma_{n_k}, \tau_{n_k}) \to(\gamma, \tau) \in\PATH$ in
the metric $\varrho$. Let $t_0(T_n)$ and $t_0(T)$ be the starting
times of the tubes $T_n$ and $T$, respectively. Their ending times
$t_1(T_n)$ and $t_1(T)$ are defined similarly. Note that $ t_i(T_{n_k})
\to t_i(T)$, for $ i=0,1$.

Since $\tau_{n_k} \le t_0(T_{n_k})$, we have $\tau\le t_0(T)$. We
first need to show that $(\gamma(s), s) \in[T]$ for each $s \in
(t_0(T), t_1(T))$. Fix $s \in(t_0(T), t_1(T))$. Since we can find
$s_{n_k} \in[t_0(T_{n_k}), t_1(T_{n_k})]$ such that $s_{n_k} \to s$.
Then $(\gamma_{n_k}(s_{n_k}), s_{n_k}) \in[T_{n_k}]$. By uniform
convergence of $\gamma_{n_k}$ to $\gamma$ on the compact time
intervals, we have $\gamma_{n_k}(s_{n_k}) \to\gamma(s)$, and hence
$(\gamma(s), s) \in[T]$. Since $(\gamma_{n_k}(t_i(T_{n_k})),
t_i(T_{n_k})) \to(\gamma(t_i(T)), t_i(T))$ for each $i$, we can
easily see that $(\gamma(t_i), t_i) \in\partial_i T$ for $i=1,2$.
This completes the proof.
\end{pf}

\begin{remark}\label{Rcompactcollection} The conclusion of the above
lemma holds under the following slightly weaker assumption which turns
out to be convenient in practice.
Let $T$ be a tube and for $(\gamma, t_0) \in\Pi$ set $ \tau(\gamma
, T):= \inf\{t\ge t_0\dvtx  \gamma(t) \in\partial_0 T \}$ (with $\inf
\varnothing= +\infty$, as usual), and set $\gamma_T(t) = \gamma(t)$
for $t\ge\tau(\gamma, T)$. If $\PATH\subseteq\Pi$ set
$\PATH_T = \{ (\gamma_T, \tau(
\gamma, T))\dvtx  \tau(\gamma,T)<\infty\}$.
Then the conclusion of Lemma~\ref{compactimpliescont} holds as soon
as $\PATH_T$ is compact in $\Pi$ for every tube $T \in\TUBE$.
In words, if the collection of paths $\PATH$, restricted to any
particular tube $T$, is compact, then this induces an element $\Cr
(\PATH) \in\HH$.
\end{remark}


\subsection{Countable coalescing process} \label{subsetcoalrule}
Fix $d \ge1$ and a countable ordered set of $\R^d$-valued continuous
paths $ (\gamma_j, s_j  )_{ j \ge1}$ in $\Pi$. From this
countable set of {\em free} paths, we define an ordered set of {\em
coalescing} paths $ (\gamma^c_j, s_j  )_{ j \ge1}$ in $\Pi$
inductively using the following coalescing rule. Set $\gamma^c_1 =
\gamma_1$. For $j >1$, set
\[
\tau_j = \inf \bigl\{ t \ge s_j\dvtx  \gamma_j(t)
\in\bigl\{\gamma^c_1(t), \gamma^c_2(t),
\ldots, \gamma^c_{j-1}(t) \bigr\} \bigr\}
\]
with the usual convention that $\inf\varnothing= \infty$. Take
\[
I_j = \min \bigl\{ i \in\{ 1,2, \ldots, j-1\}\dvtx
\gamma_j(\tau_j) = \gamma^c_{i}(
\tau_j) \bigr\} \qquad\mbox{if } \tau_j < \infty.
\]
For $t \ge s_j$, define
\[
\gamma^c_j(t) = \cases{ \gamma_j(t), &
\quad if $t < \tau_j$,
\cr
\gamma^c_{I_j}(t), &
\quad if $t \ge\tau_j$.}
\]
In words, if the free paths of labels $i$ and $j$ collide, they both
subsequently follow the path with the lower label.

\subsection{Characterization and convergence criterion for probability measures on~$\mathscr{H}$}

A general tube can be extremely complicated and its crossing
probability can be very hard to deal with. But in order to characterize
a probability measure on $\HH$, it is enough to know the joint
crossing probabilities of finitely many ``nice'' tubes belonging to a
class of tubes which is dense in $\TUBE$.

\begin{lemma}[(Characterization)]\label{lemcharacmeasure}
Let $\hat\TUBE$ be any dense subset of $\TUBE$. Let $\bQ_1$ and
$\bQ_2$ be two probability measures on $\HH$ such that for all $m \ge
1$ and for all $T_1, \ldots, T_m \in\hat\TUBE$,
\[
\bQ_1 ( \boxminus_{T_1} \cap\cdots\cap
\boxminus_{T_m} ) = \bQ_2 ( \boxminus_{T_1} \cap
\cdots\cap\boxminus_{T_m} ).
\]
Then $\bQ_1 = \bQ_2$.
\end{lemma}

\begin{pf}
The events of the form $\boxminus_{T_1} \cap\cdots\cap\boxminus
_{T_m}$ for some $m\ge1$, and some $T_1, \ldots, T_m \in\hat\TUBE
$, form a $\pi$-system which generates the entire Borel $\sigma
$-field by Theorem~\ref{tcompact}. By Dynkin's lemma, any two
probability measures which agree on this $\pi$-system must hence be identical.
\end{pf}

For us the main advantage of $\HH$ being compact is that it implies
that the set of all probability measures on $\HH$ is also compact
under the topology of weak convergence, hence any sequence of
probability measures on $\HH$ automatically has a subsequential limit.
This greatly reduces the amount to work necessary to establish weak
convergence for probability measures on $\HH$. The following
proposition gives a useful criterion for the weak convergence of a
sequence of probability measures $(\bQ_\eta)_{\eta>0}$ toward a
limiting measure $\bQ$ in terms of the joint crossing probabilities of
certain family of tubes. But before that, we state the following definition.

A subset of tubes $\hat\TUBE\subseteq\TUBE$ is called {\em
super-dense} in $\TUBE$ if:
\begin{longlist}[(P1)]
\item[(P1)] There exists a countable subset $\TUBE^0 \subseteq\TUBE$
such that $\TUBE^0$ is dense in $\TUBE$.
\item[(P2)] For each $T \in\TUBE^0$, there exists a monotone chain of
tubes $ T^\delta\in\hat\TUBE$ indexed by $\delta\in I(T, \hat
\TUBE)$, where $I (T, \hat\TUBE) \subset[0, 1]$ has a countable
complement, such that $T^0=T$ and $T^{\delta_1} < T^{\delta_2}$ if
$\delta_1 > \delta_2 \geq0$ and moreover, $d_\TUBE(T^\delta, T)
\to0$ as $ \delta\to0+$.
\end{longlist}
Clearly, a super-dense family of tubes is also dense in $\TUBE$.

\begin{proposition}[(Convergence)]\label{propconvcriteriontube}
Assume that $\hat\TUBE$ is super-dense in $\TUBE$.
Let $(\bQ_\eta)_{\eta>0}$ be a sequence of probability measures on
$\HH$. Suppose that
%
\begin{equation}
\label{alimitcrossing} q(T_1, \ldots, T_m):= \lim
_{\eta\to0} \bQ_\eta( \boxminus_{T_1} \cap\cdots
\cap\boxminus_{T_m} )
\end{equation}
exists for all $m \ge1$ and $T_1, T_2, \ldots, T_m \in\hat\TUBE$.
Then there exists a (unique) \margin{C: added unique here, even though
this was implicit by the statement itself} probability measure $\bQ$
on $\HH$, such that $\bQ_\eta\stackrel{d}{\to} \bQ$.
Furthermore there exists $\hat\TUBE' \subseteq\hat\TUBE$ which is
again super-dense in $\TUBE$ and such that $\bQ( \boxminus_{T_1}
\cap\cdots\cap\boxminus_{T_m}) = q(T_1, \ldots, T_m) $ for all
$T_1, \ldots, T_m \in\hat\TUBE'$.
\end{proposition}


Before proving the proposition, we will first note a result about the
(topological) boundary of the event that a tube is traversed. If $A
\subset\HH$, let $\partial A$ denote the boundary of $A$ for the
topology of $\HH$, that is, $\partial A = \bar A \setminus A^\circ$.

\begin{lemma} \label{lemboundarytube}
For $T \in\TUBE$, the boundary of the closed set $\boxminus_T$
satisfies \margin{A: changed the notation $A^\complement$ to $\neg A$
to keep it consistent with lemma A.2 }
\[
\partial\boxminus_T \subseteq \biggl(\bigcap
_{ T < T'} \neg\, \boxminus_{T'} \biggr) \cap
\boxminus_T.
\]
\end{lemma}

\begin{pf} 
Fix $T'$ such that $T < T'$. Find an open neighborhood $U$ of $T'$ such
that $T < T''$ for all $T'' \in U$. Then $\boxminus^U$ is an open
subset of $\boxminus_T$ and hence $\boxminus_{T'} \subseteq\boxminus
^U \subseteq\boxminus_T^\circ$. Hence, since $\boxminus_T$ is closed,
\[
\partial\boxminus_T = \overline{\boxminus_T} \setminus
\boxminus _T^\circ\subseteq\boxminus_T \cap\neg
\boxminus_{T'}.
\]
The lemma follows since the above inclusion holds for all $T< T'$.
\end{pf}


\begin{pf*}{Proof of Proposition~\ref{propconvcriteriontube}}
By compactness of $\HH$, it suffices to establish the unique
subsequential limit is $\bQ$. Let $\bQ'$ be another subsequential
limit of $(\bQ_\eta)$. Note that by the portmanteau theorem,
\[
\bQ(\boxminus_{T_1} \cap\cdots\cap\boxminus_{T_m} )= \lim
_{\eta
\to0+}\bQ_\eta( \boxminus_{T_1} \cap\cdots
\cap\boxminus_{T_m} ) = \bQ' (\boxminus_{T_1}
\cap\cdots\cap\boxminus_{T_m} )
\]
for all $m \ge1$ and $T_1, \ldots, T_m \in\hat\TUBE$ such that
$\bQ(\partial\boxminus_{T_i}) = 0$ and $\bQ'(\partial\boxminus
_{T_i}) = 0$ for each $1 \le i \le m$.

Let $\TUBE^0 \subseteq\TUBE$ be a countable subset satisfying the
properties \textup{(P1)} and \textup{(P2)} in the definition of the super-dense class of tubes.
By Lemma~\ref{lemboundarytube}, for any two tubes $T_1 < T_2$,
\[
\partial \boxminus_{T_1} \cap\,\partial\boxminus_{T_2} =
\varnothing.
\]
Let $T \in\TUBE^0$ and let $(T^\delta)_{ \delta\in I(T, \hat\TUBE
)} \subseteq\hat\TUBE$ satisfy the property \textup{(P2)} of super denseness.
Therefore, by $\sigma$-additivity of probability measures,
$\bQ(\partial\boxminus_{T^\delta}) >0$ or $\bQ'(\partial\boxminus
_{T^\delta}) >0$ can be true for only countably many $\delta\in I(T,
\hat\TUBE)$.
This implies that
\[
\hat\TUBE': = \bigl\{ T\in\hat\TUBE\dvtx  \bQ(\partial\boxminus
_{T}) = 0, \bQ'(\partial\boxminus_{T}) = 0
\bigr\}
\]
is again super-dense in $\TUBE$.
On the other hand, for all $m \ge1$ and for all $T_1, \ldots, T_m \in
\TUBE'$,
\[
\bQ(\boxminus_{T_1} \cap\cdots\cap\boxminus_{T_m} ) =
\bQ' (\boxminus_{T_1} \cap\cdots\cap\boxminus_{T_m}
) = q(T_1, \ldots, T_m).
\]
Hence, by Lemma~\ref{lemcharacmeasure}, $\bQ' = \bQ$, as desired.
\end{pf*}

\subsubsection{Example: A super-dense family of ``nice'' tubes in \texorpdfstring{$\mathbb{R}^d\times\mathbb{R}$}{{R}dx{R}}} \label{dktubes}

\begin{definition}\label{Dnicetubes}
Let $\boxT$ be the family of all tubes $T$ such that:
\begin{longlist}[1.]
\item[1.] The set $[T]$ can be tiled by a {\em finite} number of boxes of
the form $[a_1, b_1] \times\cdots\times[a_{d+1}, b_{d+1}]$ with
$a_i<b_i$ for $1\le i \le d+1$, and

\item[2.] $\partial_0 T = [T] \cap(\R^d \times\{t_0\})$ and $\partial
_1 T = [T] \cap(\R^d \times\{t_1\})$, where $t_0$ and $t_1$ are the
start and the end time of $T$.
\end{longlist}
\end{definition}

Let us now check that $\boxT$ is super-dense. Define $\TUBE^0$ same
as above the definition of $\boxT$ with an added restriction that
$a_i$ and $b_i$ appearing in item 1 \margin{Changed (D.a) ?? into item
1} are all rationals.
The family of tubes $ \TUBE^0$ is clearly countable and is dense in
the space of tubes $\TUBE$ endowed with the above metric $d_\TUBE$.
It remains to check property \textup{(P2)}. Let us fix a tube $T$ in $ \TUBE
^0$. Suppose that
\[
\partial_0 T = \prod_{i=1}^d
[a_i, b_i] \times\{ t_-\}\quad\mbox {and}\quad
\partial_1 T = \prod_{i=1}^d
[c_i, d_i] \times\{ t_+\}.
\]
By the definition of the class of tubes $\boxT$, we can find $t_e>0$
such that
\begin{eqnarray*}
[T] \cap\bigl(\R^d \times[t_-, t_- + t_e]\bigr) &=& \prod
_{i=1}^d [a_i,
b_i] \times[t_-, t_- + t_e],
\\
{[T]} \cap\bigl(
\R^d \times[t_+ - t_e, t_+]\bigr) &=& \prod
_{i=1}^d [c_i, d_i]
\times[t_+ - t_e, t_+].
\end{eqnarray*}

For $\delta< t_e$, define $T^\delta$ to be the tube such that
$[T^\delta] = [T]^\delta\cap (\R^d \times[t_- +\delta, t_+
-\delta] )$, $\partial_0 T^\delta= \prod_{i=1}^d [a_i -\delta,
b_i+\delta] \times\{ t_- +\delta\}$ and\vspace*{1pt} $\partial_1 T^\delta=
\prod_{i=1}^d [c_i-\delta, d_i+\delta] \times\{ t^+ -\delta\}$,
where for any set $S \subseteq\R^{d+1}$, the enlargement $S^\delta$
is defined to be the closed set consisting of all points in $\R^{d+1}$
whose $L^\infty$ distance from $S$ is at most $\delta$.
It is easy to check that $T^0 = T$, $T^{\delta_2} < T^{\delta_1}$ if
$0 \le\delta_1 < \delta_2 < t_e$ and $d_\TUBE(T^\delta, T) \to0$
as $\delta\to0+$. This shows the property \textup{(P2)} holds.

\section{Arratia's flow in $\mathscr{H}$}

\label{SBW}

\subsection{Existence}\label{SSArratiaflow}
In this section, we restrict our attention to Brownian motion on $\R$
and define a unique measure on $\HH$ which represents the coalescing
Brownian flow on $\R$ in the tube topology (for $d=1$).
We start by introducing some notation. Let $z_1 = (x_1, t_1),z_2 =
(x_2, t_2), \ldots$ be a sequence of space--time points in $\R^2$ and
let $\cD= \{z_1, z_2, \ldots\}$ be a countable ordered set. We assume
that $\cD$ is dense in $\R^2$. Let $(W_j(t))_{t \ge t_j}$ be
independent Brownian motions starting, respectively, from the space--time
points $(z_j)_{j\ge1}$, that is, $W_j(t_j) = x_j$, defined on a common
probability space $(\Omega, \mathcal{F}, \prob)$.
Using the coalescing rule as described in Section~\ref{subsetcoalrule}, for each $n\ge1$ this defines
a collection of $n$ coalescing paths denoted by $W^c_1, W^c_2, \ldots,
W^c_n$ starting from $z_1, z_2,\ldots, z_n$, respectively. This
collection of coalescing paths, being finite, is of course compact in
$(\Pi, \varrho)$, and hence induces, by Lemma~\ref
{compactimpliescont}, a random element $ \calW_n \in\HH$ defined by
%
\begin{equation}
\label{eWn} \calW_n:= \calW(z_1, \ldots,
z_n): = \Cr\bigl(\bigl\{ W^c_1,
W^c_2, \ldots, W^c_n\bigr\}
\bigr).
\end{equation}
We now state the main theorem of the section.

\begin{theorem}[(Coalescing Brownian flow on $\R$)] \label{tBW}
The random variables $\cW_n$ converge in distribution as $n\to\infty
$ to a random variable $\cW_\infty$ in $\HH$\margin{C: added $\P
_\infty$ in the statement since it did not appear before}, whose law
$\P_\infty$ does not depend on the dense countable set $\cD$
(including its order).
\end{theorem}

\begin{definition}
\label{dBW} A random variable on $\HH$ with law $\P_\infty$ is
called a {\em coalescing Brownian flow} on $\R$.
\end{definition}

\begin{pf}
By construction, $\calW_1 \subseteq\calW_2 \subseteq\cdots$ almost
surely. That is, $\cW_k$ is a nondecreasing sequence. Hence, by Lemma
\ref{lmonotone}, it follows that the $\lim_{k \to\infty} \cW_k$
exists almost surely in $\HH$ and is equal to $\overline{\bigcup
_{k=1}^\infty\calW_k}$. We call $\cW(\D)$ the limiting element of
$\HH$ and denote by $\P_\infty^\D$ its law on $\HH$ which, at this
stage, might depend on the ordered set $\D$. We aim to show that the
law $\P_\infty^\D$ does not in fact depend on $\cD$.


%


\begin{lemma} \label{Dindep}
Let $\mathcal D = ( z_j )_{ j \ge1} $ and $\mathcal D' = ( z'_j )_{ j
\ge1} $ be two countable dense ordered sets of $\mathbb R^2$.
Fix $m \ge1$ and tubes $T_1, \ldots, T_m \in\boxT$. Let
\begin{eqnarray*}
p(T_1, \ldots, T_m) &=& \lim_{n \to\infty}
\prob\bigl(\cW(z_1, \ldots, z_n) \in\boxminus_{T_1}
\cap\cdots\cap\boxminus_{T_m}\bigr)\quad\mbox {and}
\\
p'(T_1, \ldots, T_m) &=& \lim
_{n \to\infty} \prob\bigl(\cW\bigl(z'_1,
\ldots, z'_n\bigr) \in\boxminus_{T_1}\cap\cdots
\cap\boxminus_{T_m}\bigr).
\end{eqnarray*}
Then
$
p(T_1, \ldots, T_m) = p'(T_1, \ldots, T_m)$.
\end{lemma}

\begin{pf}
The limits $p(T_1, \ldots, T_m)$ and $p'(T_1, \ldots, T_m)$ exist due
to monotonicity.
Let $z_j = (x_j, t_j)$ and $z'_j = (x'_j, t_j')$.
We\vspace*{1pt} may suppose without loss of generality that the points
$z_1, z_2, \ldots$ (resp., $z_1', z_2', \ldots$) are distinct. For
$1 \le i \le m$, let the lower face of the tube $T_i$ is given by
$\partial_0 T_i = [a_i, b_i] \times\{ u_i\}$, $u_i$ being the start
time of $T_i$. We denote by $A$ the event $ \boxminus_{T_1} \cap
\cdots\cap\boxminus_{T_m}$.


Note that a Brownian particle starting from a point in the boundary of
$[a_i, b_i]$ at time $u_i$ will immediately escape the interval $[a_i,
b_i]$ almost surely and the tube $T_i$ can not be crossed by it. So, as
far as the event $A$ is concerned we can assume that none of the points
in $z_j$ or $z'_j$ lies in $\{ (a_i, u_i), (b_i, u_i)\dvtx  1 \le i \le n\}
$, the boundary of the lower face of any tube.

For $x, t \in\mathbb R$ and $\eps>0$, define two rectangles $R(x, t;
\eps) = [x-\eps, x+ \eps] \times[t - \eps, t+\eps]$ and $R^+(x,
t; \eps) = [x-\eps, x+ \eps] \times[t, t+\eps]$.
Let $W$ be Brownian motion starting from the space--time point $(x, t)$
and $W'$ be another independent Brownian motion starting from the space--time point $(x', t')$. Let $\tau\ge\max(t, t')$ be the first hitting
time of $W$ and $W'$. Then for any fixed $\eps>0$, the probability of
the event
%
\begin{eqnarray}\label{eqimmediatehit}
&& \bigl\{ \tau< t+\eps\mbox{ and }\bigl( W(s), s \bigr) \in R(x, t;
\eps)\ \forall t \le s \le\tau\mbox{ and}
\nonumber\\[-8pt]\\[-8pt]\nonumber
&&\hspace*{81pt} \bigl( W'(s), s
\bigr) \in R(x, t; \eps)\ \forall t' \le s \le\tau \bigr\}
\end{eqnarray}
\margin{C: added the condition $\tau< t+\eps$ in the event}
converges to $1$ as $(x', t') \to(x, t)$. This follows from the right
continuity of Brownian paths and the fact that $\tau\searrow t$ almost
surely as $(x', t') \to(x, t)$.

Fix $n \in\mathbb N$ and $\delta>0$.
We can choose $\eps_0>0$ sufficiently small such that it satisfies:
\begin{longlist}[(iii)]
\item[(i)] $R(x_j, t_j, \eps_0), 1 \le j \le n$ are pairwise disjoint.
\item[(ii)] For each $1 \le j \le n$ and each $1 \le i \le m$, if $t_j
< u_i$ then $t_j + \eps_0 < u_i$.
\item[(iii)] For each $1 \le j \le n$ and each $1 \le i \le m$, if
$t_j = u_i$ and $x_j \in(a_i, b_i)$, then $R^+(x_j, t_j, \eps_0)
\subseteq[T_i]$.
\end{longlist}

Given any $\delta>0$ and some $\eps> 0$ that will be specified in a
moment, find
$y_1' = (v_1', s_1'), \ldots, y_n' = (v_n', s_n') \in\D'$
such that $s_j' \le t_j$ and $\llvert  x_j - v_j'\rrvert   + \llvert   t_j - s'_j\rrvert   \le\eps$
for $1 \le j \le n$. This is possible since $\D'$ is dense in $\mathbb R^2$.
By (\ref{eqimmediatehit}), there is $\eps_0$ such that
if we choose $\eps< \eps_0$ sufficiently small, then for each $1 \le
j \le n$, two independent Brownian motions starting from the space--time
points $z_j$ and $y'_j$
collide before the graphs of their trajectories leave the rectangle
$R(x_j, t_j, \eps_0)$ with probability at least $1-\tfrac{\delta}{
n}$. \margin{C: changed at most by at least $1-$.} Obviously,
%
\begin{equation}
\label{eqBWdensetset1} \prob\bigl( \calW(z_1, \ldots, z_n) \in A
\bigr) \le\prob\bigl( \calW\bigl( y_1', \ldots,
y_n', z_1, \ldots, z_n\bigr)
\in A \bigr),
\end{equation}
by invariance under reordering (using the strong Markov property of
Brownian motion).
Thus, by our choice of $\eps$ and a simple union bound, we obtain
%
\begin{eqnarray}
\label{eqBWdensetset3}
&& \prob\bigl( \calW\bigl( y_1', \ldots,
y_n', z_1, \ldots, z_n\bigr)
\in A \bigr)\nonumber
\\
&&\qquad \le \prob\bigl( \calW\bigl( y_1', \ldots,
y_n'\bigr) \in A \bigr) +\delta
\\
&&\qquad \le p'(T_1, \ldots, T_m) + \delta.\nonumber
\end{eqnarray}
Combining (\ref{eqBWdensetset1}) and (\ref{eqBWdensetset3}), we obtain that, for each $n \ge1$
\[
\prob\bigl(\cW(z_1, \ldots, z_n) \in A\bigr) \le
p'(T_1, \ldots, T_m) +\delta. %
\]
Taking limit as $n \to\infty$ and noting that $\delta>0$ is
arbitrary, we have $p(T_1, \ldots,\break T_m) \le p'(T_1, \ldots,  T_m)$.
Interchanging the role of $\D$ and $\D'$ we deduce the equality. This
completes the proof of the lemma.
\end{pf}

%

Now, let us show how Lemma~\ref{Dindep} implies Theorem~\ref{tBW}.
Take any countable dense set $\D' = (z_j')_{j \ge1}$. By Lemma~\ref
{Dindep}, for fixed tubes $T_1, \ldots, T_m \in\boxT$,
\begin{eqnarray*}
p(T_1, \ldots, T_m ) &=& \lim_{n \to\infty}
\prob\bigl(\calW(z_1, \ldots, z_n) \in A \bigr)= \lim
_{n \to\infty} \prob \bigl( \calW\bigl(z'_1,
\ldots, z'_n\bigr) \in A \bigr)
\\
&=& p'(T_1,
\ldots, T_m).
\end{eqnarray*}
Since $\boxT$ is super-dense in $\TUBE$, by Proposition~\ref
{propconvcriteriontube}, we know that $p(T_1, \ldots,\break  T_m) = \P
^{\cD}_\infty(A)$ and $p'(T_1, \ldots, T_m) = \P^{\cD'}_\infty
(A)$ at least when $T_1, \ldots, T_m \in\TUBE'$, where $ \TUBE'$ is
super-dense in $\TUBE$. Thus, we conclude $\P^{\cD}_\infty(A) = \P
^{\cD'}_\infty(A)$ if $T_1, \ldots, \break T_m \in\TUBE'$. By Lemma~\ref
{lemcharacmeasure}, this shows $\P^{\cD}_\infty= \P^{\cD
'}_\infty$ and completes the proof of Theorem~\ref{tBW}.
\end{pf}

\subsection{Characterization of Arratia's flow}\label{Scharac}

Let $\P_\infty$ be the law of $\cW_\infty$ on $\HH$, and for any
tubes $T_1, \ldots, T_m \in\TUBE$ (and hence in $\boxT$), let
\[
p(T_1, \ldots, T_m) = \lim_{n \to\infty}
\prob\bigl(\cW(z_1, \ldots, z_n) \in\boxminus_{T_1}
\cap\cdots\cap\boxminus_{T_m}\bigr).
\]
It follows from Theorem~\ref{tBW} and Proposition~\ref
{propconvcriteriontube} that there exists a super-dense family of
tubes $ \boxT'$ such that for all $T_1, \ldots, T_m \in\boxT'$,
%
\begin{equation}
\label{eqBW} p(T_1, \ldots, T_m) = \P_\infty(
\boxminus_{T_1}\cap\cdots\cap \boxminus_{T_m}).
\end{equation}
This characterises uniquely $\P_\infty$, though in practice a
drawback of this conclusion is that we do not know what $ \boxT'$ is.
However, the following result shows that this conclusion remains valid
for all $T_1, \ldots, T_m \in\boxT$. We stress, however, that this
result is not needed for the rest of the paper, so this section may be
skipped by a reader who is only interested in the invariance principle
(Theorem~\ref{tscalinglimit}).

\begin{theorem}\label{tBW2}
$\P_\infty$ is the unique probability distribution on $\HH$ such
that for all $m \ge1$ and for all fixed tubes $T_1, \ldots, T_m \in
\boxT$,
%
\begin{eqnarray}\label{eqBWcharac}
&& \P_\infty(\boxminus_{T_1}\cap\cdots\cap
\boxminus_{T_m} )
\nonumber\\[-8pt]\\[-8pt]\nonumber
&&\qquad  = \sup_{{n \ge1; z_1, \ldots, z_n \in\R^2}} \prob \bigl(
\calW(z_1, \ldots, z_n) \in\boxminus_{T_1}\cap
\cdots\cap\boxminus_{T_m} \bigr).
\end{eqnarray}
\end{theorem}

\begin{pf}
Note that the supremum in the right-hand side of (\ref{eqBWcharac})
is simply $p(T_1, \ldots, T_m)$. The lower bound is easy: indeed,
since $A = \boxminus_{T_1} \cap\cdots\cap\boxminus_{T_m}$ is
closed in $\HH$, and since $\cW(z_1, \ldots, z_n) \to\cW_\infty$
in distribution on $\HH$,
\[
\P_\infty(A) \ge\lim_{n\to\infty} \prob\bigl(
\cW(z_1, \ldots, z_n) \in A\bigr). %
\]
Taking a supremum over $n \ge1$ and $z_1, \ldots, z_n$ shows that
\[
\P_\infty(A ) \ge\sup_{{n \ge1; z_1, \ldots, z_n \in\R^2}} \prob\bigl(
\calW(z_1, \ldots, z_n) \in A\bigr). %
\]

We now turn to proof of the upper bound in (\ref{eqBWcharac}).
%
Recall that by Lemma~\ref{lmonotone}, $\P^\D_\infty(A) = \prob
(T_1, \ldots, T_m \in\overline{\bigcup_k \calW_k} )$, where $A =
\boxminus_{T_1} \cap\cdots\cap\boxminus_{T_m}$.
We need to prove that for any finite number of tubes $T_1, T_2, \ldots
, T_m \in\boxT$,
\[
\prob \Bigl(T_1, \ldots, T_m \in\overline{\bigcup_k
\calW_k} \Bigr) = \prob \Bigl(T_1, \ldots, T_m \in
\bigcup_{k} \calW_k \Bigr),
\]
as the right-hand side of the above equation
is the increasing limit (as $k \to\infty$) of
\[
\prob \bigl( \calW(z_1, \ldots, z_k) \in
\boxminus_{T_1}\cap\cdots \cap\boxminus_{T_m} \bigr).
\]
Note that it is enough to prove the above equality for $m=1$ since
\begin{eqnarray*}
&& \prob \Bigl(T_1, \ldots, T_m \in\overline{\bigcup_k
\calW_k} \Bigr) - \prob \Bigl(T_1, \ldots, T_m \in
\bigcup_{k} \calW_k \Bigr)
\\
&&\qquad \le\sum_{i=1}^m
\Bigl(\prob (T_i \in\overline{\bigcup_k
\calW_k} ) - \prob (T_i \in\calW_k )
\Bigr).
\end{eqnarray*}
So, we fix a tube $T \in\boxT$. We are going to show that
%
\begin{equation}
\label{eqonetubelimitcrossing} \prob \Bigl(T \in\bigcup_{k} \calW_k \Bigr) = \prob \Bigl(T
\in \overline{\bigcup_k \calW_k} \Bigr).
\end{equation}
Now $T \in\overline{\bigcup_k \calW_k}$ means that there exists an
increasing sequence $k_n$ of integers and tubes $T_n \in\calW_{k_n}$
such that $T_n \to T$. For any tube $T' < T$, we can always find an
open neighborhood $ U$ of $T$ such that $T' < T''$ for all $T'' \in U$.
Consequently, $T' < T_n$ for large enough $n$ and thus $ T' \in\calW
_{k_n} \subseteq\bigcup_k \calW_k$. This implies that for any tube $T'
< T$, we have
\[
\prob \Bigl(T \in\overline{\bigcup_k \calW_k} \Bigr) \le\prob
\Bigl(T' \in\bigcup_k \calW_k \Bigr).
\]
Our goal is to show that for any $\eps>0$, there exists a tube $T'
\in\TUBE$, such that $T' < T$ and
%
\begin{equation}
\label{eqonetubeinequality} \prob \Bigl(T' \in\bigcup_k
\calW_k \Bigr) \le\prob \Bigl(T \in\bigcup _k
\calW_k \Bigr) + \eps,
\end{equation}
which immediately implies (\ref{eqonetubelimitcrossing}).
Let $\partial_0 T = [a, b] \times\{ t_0\}$ and $\partial_1 T = [c,
d] \times\{ t^*\}$.
\margin{C: why $t^*$ ?? I was confused by this notation, can I changed
to $t_1$ ?} By the definition of the class of tubes $\boxT$, we can
find $t_e>0$ such that
$[T] \cap(\R\times[t_0, t_0 + t_e]) = [a, b] \times[t_0, t_0 +
t_e]$ and $[T] \cap(\R\times[t^* - t_e, t^*]) = [c, d] \times[t^* -
t_e, t^*]$.

For $\delta< t_e$, recall the tube $T^\delta\in\boxT$ introduced in
Section~\ref{dktubes}: this is the tube such that $[T^\delta] =
[T]^\delta\cap (\R\times[t_0 +\delta, t^* -\delta] )$,
$\partial_0 T^\delta= [a -\delta, b+\delta] \times\{ t_0 +\delta\}
$ and $\partial_1 T^\delta= [c-\delta, d+\delta] \times\{ t^*
-\delta\}$.
It follows from the definition of the class of tubes $\hat\TUBE$ that
$T^\delta< T$ and $d_\TUBE(T^\delta, T) \to0$ as $\delta\to0+$.

We are going to prove (\ref{eqonetubeinequality}) taking $T' =
T^\delta$ for some small enough $\delta$, depending only on $\eps$
and $T \in\boxT$.
We are going to argue this by a (rather long) series of simple observations.

For $s<t$ and an interval $I$, let $A_I(s,t)$ denote the locations of
particles at time $t$ whose trajectories started before or at time $s$
and stayed in $I$ throughout $[s,t]$. The proof of Lemma~\ref{Dindep}
also shows that the law of $A_I(s,t)$ does not depend on~$\cD$.
We first need a well-known and simple lemma which states that the
coalescing Brownian paths ``come down from infinity''.

\begin{lemma}
\label{Lcdi} If $I $ is bounded then $\llvert  A_I(s,t) \rrvert  <\infty$ almost
surely for $s< t$. Its law depends only on $t-s$ and $I$, and is
continuous (say in total variation) in both $s$ and $t$ if $s<t$.
\end{lemma}

\begin{pf}
The first statement is well known and follows, for instance, from
Arratia's work \cite{arratia79}. The continuity of the law of
$\llvert  A_I(s,t)\rrvert  $, in total variation, is trivial to verify as $s$ is fixed
and $t$ varies, since a.s. $\lim_{\eta\to0} \llvert  A_I(s, t +\eta)\rrvert   =
\llvert  A_I(s,t)\rrvert  $. Hence, continuity follows in both $s$ and $t$ provided
that $s<t$.
\end{pf}

Now let $\alpha= \frac{t_e}{4}$ be fixed. Set $I = [a,b]$ and
$I^\delta= [a-\delta, b+\delta]$, where $\delta$ will be chosen
sufficiently small.
Let $t_1 = t_0 + \delta$ and $t_2 = t_0 + \delta+ \alpha$. Consider
the sets $A_0 = A_I(t_0, t_2)$ and $A_1 = A_{I^\delta}(t_1, t_2)$.

\begin{lemma}
\label{LA2}
For $\eps>0$, we can choose $\delta_0 >0$ so that $A_0$ and $A_1$
agree with probability greater than $1-\eps$ for all $\delta< \delta_0$.
\end{lemma}

\begin{pf}
Note that $A_0 \subset A_1$, and hence it suffices to show that $\prob
( \llvert  A_0\rrvert   = \llvert  A_1\rrvert  ) \ge1- \eps$.
It is a straightforward consequence of Lemma~\ref{Lcdi} that
$\llvert  A_I(t_1, t_2)\rrvert  $ and $\llvert  A_I(t_0, t_2)\rrvert  $ agree with probability greater
than $1-\eps/2$ for $\delta$ sufficiently small. Using scale and
translation invariance of Brownian motion and the same argument, we see
also that $\llvert  A_I(t_1, t_2)\rrvert  $ and $\llvert  A_{I^\delta}(t_1, t_2)\rrvert  $ agree with
probability greater than $1-\eps/2$. Hence, the result follows.
\end{pf}
Now for $\kappa>0$, define $J= [a+\kappa, b-\kappa]$.

\begin{lemma}\label{LA1}
For all $\eps>0$, we can choose $\kappa>0$ and $\delta_0>0$ so that,
uniformly in $ \delta< \delta_0$, the tube $T^\delta$ is crossed if
and only if it is crossed by a path touching $(A_1 \cap J ) \times\{
t_2\}$, except with probability at most $\eps$.
\end{lemma}

\begin{pf} Obviously, if $T^\delta$ is crossed by a particle, then
that particle has to touch $A_1 \times\{t_2\}$. So, it suffices to
show that with probability at least $1 -\eps$, none of the particles
starting from $(A_1 \cap J^c ) \times\{ t_2\}$ will stay inside
$[T^\delta]$ up to time $t_0+ \tfrac{t_e}{2}$.

Clearly, $I^\delta\setminus J \subset[a - \kappa, a+\kappa] \cup[b
- \kappa, b+\kappa]$ if $\delta< \kappa$. Now, choose $n$ large
enough that for all $\delta< \tfrac{\alpha}{2}$, $\llvert  A_1 \rrvert   \le
\llvert  A_{[a-1, b+1]}(t_0+ \tfrac{\alpha}{2}, t_0+ \alpha)\rrvert   \le n$ with
probability greater than $1- \tfrac{\eps}{3}$. Observe that we can
take $\kappa>0$ small enough such that any coalescing Brownian path in
$[a - \kappa, a+\kappa]$ at time $t_2$ will hit the line $x = a -
\kappa$ by time $t_2 + \tfrac{t_e}{4} $ with probability at least $1
- \tfrac{\eps}{3n}$, uniformly in $\delta$. Such a particle
necessarily leaves $T$ and $T^\delta$, if $\delta< \kappa$.
Likewise, any coalescing Brownian path in $[b-\kappa, b+\kappa]$ at
$t_2$ will also hit the line $y=b +\kappa$ by time $ t_2 + \tfrac
{t_e}{4}$ with probability at least $1 - \tfrac{\eps}{3n}$, uniformly
in $\delta$.
Summing over all particles at time $t_2$,
we see that with probability greater than $1-\eps$, any particle at
time $t_2 $ located within $(I^\delta\setminus J) \cap A_1$ can cross
$T^\delta$.
Lemma~\ref{LA1} follows with $\delta_0 = \min(\kappa, \tfrac
{\alpha}{2})$.
\end{pf}

Combining Lemmas~\ref{LA2} and~\ref{LA1}, we deduce that for
$\delta< \delta_0$ (with $\delta_0$ as in Lemmas~\ref{LA2} and
\ref{LA1}, depending only on $\eps$ and $T$),
\[
\prob\Bigl(T^\delta\in\bigcup_k \calW_k \Bigr) -
\prob\Bigl(T \in\bigcup_k \calW _k \Bigr) \le2 \eps+ p,
\]
where $p$ is the probability that one of the particles passing through
$(x, t_2 )$ for some $x \in A_1 \cap J$ will stay within $[T^\delta]$
until the time $t^* - \delta$ but the trajectory of that particle will
leave the tube $T$ at sometime between $t_2$ and $t^*$. Recall also
that, by Lemma~\ref{Lcdi}, we can find a large $n$ (depending only on
$\eps$ and $T$) such that\vadjust{\goodbreak} $\llvert  A_1\rrvert   \le n$ with probability at least $1
- \eps$. Hence, using the Markov property of Brownian motion, we
complete the proof of (\ref{eqonetubeinequality}) using the
following lemma.

\begin{lemma} \label{lemcorner}
Let $B_{\delta, x}$ be the event that a Brownian motion starting from
$x$ at time $t_2$ will stay within $T^\delta$ until $t^* - \delta$
but the trajectory of that Brownian motion leaves the tube $T$ at
sometime between $t_2$ and $t^*$. Then
%
\begin{equation}
\label{eqBdeltax} \limsup_{\delta\to0+} \sup_{ x \in[ a + \kappa, b- \kappa]}
\prob(B_{\delta, x}) = 0.
\end{equation}
\end{lemma}

\begin{pf}
Note that if $T \in\boxT$, then the boundary of the set $[T]$ can be
expressed as the union of finitely many vertical and horizontal line segments.
We call a point $(y, t)$ a {\em corner point} of $T$ if $(y, t)$ lies
at the point of intersection of a vertical and a horizontal line
segment on the boundary of $[T]$ as described in Figure~\ref{fig3} above.
Let $S$ be the set of all corner points $(y, t)$ of $T$ such that $t_2<
t \le t^*$. Let $W^x$ be Brownian motion starting from $x$ at time $t_2$.

%
\begin{figure}

\includegraphics{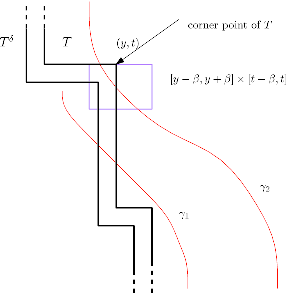}

\caption{The path $\gamma_1$ leaves both $T$ and $T^\delta$ whereas
the path $\gamma_2$ exits from $T$ near a corner point and continues
to stay inside $T^\delta$.}\label{fcornerpoint}\label{fig3}
\end{figure}

To estimate $\prob(B_{\delta, x})$, we consider the following event:
\[
W^x(s) \in[y - \beta, y + \beta]\qquad\mbox{for some } s \in[ t -
\beta, t]\mbox{ and }(y, t) \in S.
\]
%
Since the set $S$ is finite, we can find $\beta>0$ small such that the
above event has probability at most $\eps$, uniformly in $x \in[ a +
\kappa, b - \kappa]$.


Note that if the above event does not occur and if $\delta\in(0,
\beta)$, then for the event $B_{\delta, x}$ to happen, $W^x$ must
exit the tube $T$ for the first time through some point $(y', t')$ on
the left or the right boundary of $[T]$ such that:

\begin{longlist}[(ii)]
\item[(i)] the vertical line segment joining the points $(y', t')$ and $(y',
t'+ \beta)$ lies on the boundary of $[T]$ and

\item[(ii)] for all time $[t', t'+\beta]$, the trajectory of $W^x$ must
continue to stay inside the lines $x = y' - \delta$ or $x = y' +
\delta$, depending on which of the two boundaries (left or right) of
$T$ it violates. By choosing $\delta>0$ small, we can make this
probability smaller than $\eps$. Thus, for small enough $\delta>0$,
\[
\sup_{x \in[ a + \kappa, b - \kappa]} \prob(B_{\delta, x}) \le 2\eps,
\]
which proves (\ref{eqBdeltax}). This establishes (\ref
{eqonetubeinequality}), and thus completes the proof of the upper bound.\quad\qed
\end{longlist}\noqed
\end{pf}

In turn, this completes the proof of Theorem~\ref{tBW2}.
\end{pf}

\section{Invariance principle for coalescing random walks}\label{Sscalinglimit}

Consider a system of independent coalescing random walks started from
every space--time point $(x,t)$ on $\Z\times\Z$. We assume that the
step distribution $\xi$ satisfies
%
\begin{equation}
\label{axi} \E[\xi] =0,\qquad \E\bigl[\xi^2\bigr] = \sigma^2
\quad\mbox{and}\quad\xi \mbox{ is aperiodic}.
\end{equation}
Under diffusive scaling, this gives rise to a collection of continuous
paths in $(\Pi, \varrho)$ obtained by interpolating linearly the
paths of the coalescing random walk in the rescaled lattice $ \LL_\eta
: = \sigma^{-1} \eta\Z\times\eta^2 \Z$, which we will denote by
$\Gamma^\eta$.
Note that in~$\Gamma^\eta$, two paths can cross over each other
several times before they finally merge at some point in $\LL_\eta$.
Note that if $\calW^\eta:= \Cr(\Gamma^\eta) $ then $\cW^\eta\in
\HH$ by Remark~\ref{Rcompactcollection}: indeed, in the notation of
this remark, the collection of paths $(\Gamma^\eta)_T$, restricted to
any particular tube $T$, is finite, and hence compact. We call $\P
^\eta$ the law on $\HH$ of~$\cW^\eta$.

\begin{theorem}\label{tscalinglimit}
Assume (\ref{axi}). Then as $\eta\to0$,
\[
\P^\eta{\to} \P_\infty,
\]
weakly in $\HH$, where $\P_\infty$ is the law of the coalescing
Brownian flow on $\HH$, as defined in Theorem~\ref{tBW}.
\end{theorem}

\begin{remark}
As can be trivially seen from the proof, the same conclusion holds for
many variants. Here is one such example. Fix\vspace*{1pt} a step distribution $\xi$
which is nonlattice, centred and such that $\E[\xi] =0$ and $\E[\xi
^2] =\sigma^2$. Consider $\Gamma^\eta$ a system of coalescing random
walks in continuous time (jumping at rate one according to the
distribution $\xi$). The particles start from a cloud of points $(x,t)
\in\R\times\R$ distributed according to a Poisson point process
with intensity $dx\otimes dt$, and particles coalesce as soon as their
mutual distance is less than one (i.e., the path of the particle with
higher label merges with the path of the particle with lower label, in
some fixed enumeration of the Poisson cloud particles, as in
Section~\ref{SSArratiaflow}). Then applying the same diffusive
scaling, this gives rise to a law $\P^\eta$ on $\HH$ which converges
weakly to Arratia's flow $\P_\infty$ as $\eta\to0$.
\end{remark}

\begin{pf*}{Proof of Theorem \ref{tscalinglimit}}
Fix $T_1,\ldots, T_m \in\boxT$, and recall our notation from
equation~(\ref{eWn}) and Lemma~\ref{Dindep}: \margin{added these
two links here}
\[
p(T_1, \ldots, T_m) = \lim_{n \to\infty}
\prob(\cW_n \in A),\vadjust{\goodbreak} %
\]
where $A = \boxminus_{T_1}\cap\cdots\cap\boxminus_{T_m}$.
We split the proof of Theorem~\ref{tscalinglimit} into two parts, a
lower and an upper bound. The lower bound will consist in showing that
for $T_1, \ldots, T_m \in\boxT$,
%
\begin{equation}
\label{eLB} \liminf_{\eta\to0} \P_\eta(A) \ge
p(T_1, \ldots, T_m).
\end{equation}
For the upper bound, we will show that for $T_1, \ldots, T_m \in\boxT$,
%
\begin{equation}
\label{eUB} \limsup_{\eta\to0} \P_\eta(A) \le
\P_\infty(A).
\end{equation}
By Proposition~\ref{propconvcriteriontube}, there exists a
super-dense family $\boxT'$ such that when $T_1, \ldots, T_m \in
\boxT'$, it holds that $p(T_1, \ldots, T_m) = \P_\infty(A)$.
Consequently, by (\ref{eLB}) and (\ref{eUB}), if $T_1, \ldots, T_m
\in\boxT'$, then
\[
\lim_{\eta\to0} \P_\eta(A) = \P_\infty(A).
\]
Since $\boxT'$ is itself super-dense, another application of
Proposition~\ref{propconvcriteriontube} completes the proof of the theorem.

\subsection{Proof of lower bound}

Fix $n \ge1$ and let $z_1, \ldots, z_n \in\R^2$. Let $z_1^\eta,
\ldots, z_n^\eta$ be space--time points in the rescaled lattice $\LL
_\eta$ such that $z_1^\eta\to z_1, \ldots, z_n^\eta\to z_n$ as
$\eta\to0$. Let $\Gamma^\eta_n = \Gamma^\eta(z_1^\eta, \ldots,
z_n^\eta) $ be a system of $n$ independent rescaled coalescing random
walks in $\LL_\eta$ started from $z_1^\eta, \ldots, z_n^\eta$,
viewed as a random element of $(\Pi^n, \varrho^n)$ as defined in
(\ref{pathmetric}). Let $\Gamma_n = \Gamma(z_1, \ldots, z_n)$ be a
system of $n$ independent coalescing Brownian motions started from
$z_1, \ldots, z_n$, also viewed as a random variable in $(\Pi^n,
\varrho^n)$.

\begin{lemma}\label{lkpoint}
As $\eta\to0$,
\[
\Gamma_n^\eta\to\Gamma_n %
\]
in distribution on $(\Pi^n, \varrho^n)$.
\end{lemma}

The lemma says that $n$ coalescing random walks converge to $n$
coalescing Brownian motions, which is of course hardly surprising. For
a detailed proof of this fact, see \cite{nrs05}, which we will not
repeat here. But later while treating coalescing random walks on
Sierpinski gasket, we will provide a new proof of the above result that
holds in greater generality.

The proof of the lower bound (\ref{eLB}) is now easy. Fix $n \ge1$
and $z_1, \ldots, z_n \in\R^2$. Let $z_1^\eta, \ldots, z_n^\eta
\in\LL_\eta$ such that $z_i^\eta\to z_i$. Let $\cW^\eta_n = \Cr
(\Gamma^\eta(z_1^\eta, \ldots, z_n^\eta) )$ and $\cW_n = \Cr
(\Gamma(z_1, \ldots, z_n))$. By monotonicity,
\begin{eqnarray*}
\P^\eta(\boxminus_{T_1} \cap\cdots\cap\boxminus_{T_m})
&\ge& \prob\bigl( \cW_n^\eta\in\boxminus_{T_1}
\cap\cdots\cap\boxminus _{T_m}\bigr)
\\
&\to&\prob( \cW_n \in\boxminus_{T_1} \cap\cdots\cap
\boxminus_{T_m})
\end{eqnarray*}
as $\eta\to0$, where the convergence follows from Lemma~\ref{lkpoint}.
Taking a limit as $n\to\infty$, for a fixed enumeration $z_1, z_2,
\ldots$ of a dense countable set $\cD$ in $\R^2$, we obtain by
Theorem~\ref{tBW}
%
\begin{equation}
\liminf_{\eta\to0 }\P^{\eta} (\boxminus_{T_1}\cap
\cdots\cap \boxminus_{T_m}) \geq p(T_1, \ldots,
T_m),
\end{equation}
as desired.

\subsection{Uniform coming down from infinity}

Let us now prove (\ref{eUB}),
which, together with (\ref{eLB}) proves Theorem~\ref
{tscalinglimit}, as already explained. The proof of (\ref{eUB})
relies essentially on the following property, which we call \emph
{uniform coming down from infinity}. It says that, after a positive
amount of time $\delta>0$ the number of rescaled coalescing random
walks which are killed upon exiting a bounded region of space stays
finite as $\eta\to0$ (i.e., is a tight family of random variables).

\begin{proposition} \label{propcomingdownfrominf}
Let $K> 0$ be fixed. Consider coalescing random walks on $\Z\times\Z
$ with increments distributed as $\xi$ satisfying (\ref{axi}),
starting from each \mbox{$x \in[-Kn, Kn] \cap\mathbb Z$} at time $0$, and
that are killed upon leaving the interval $[ - Kn, Kn]$.
For $\delta>0$, let $U_n$ be the number of distinct coalescing random
walks at time $\delta n^2$.
Then there exists a constant $C$ independent of $\delta$ and $n$ such
that for all $k, n \in\mathbb Z_+$,
\[
\prob(U_n \ge k ) \le\frac{C}{\delta k}.
\]
\end{proposition}

We start with a simple lemma.

\begin{lemma}\label{lemmahitting}
Let $X^1, X^2$ be two independent random walks on $\Z$ with increments
distributed as $\xi$ and starting at $x, y \in\mathbb Z$ at time $0$,
respectively. Let $\tau_{x,y}$ be the integer stopping time when the
two walkers first meet. Then for all $t \in\mathbb Z_+$,
\[
\prob(\tau_{x,y} > t) \le\frac{C_0}{\sqrt t} \llvert x - y\rrvert
\]
for some constant $C_0$ independent of $t, x$ and $y$.
\end{lemma}

\begin{pf}
This bound is easy to derive and had already been used by \cite{nrs05}
(see Lemma 2.2). Assume without loss of generality that $x<y$. When
$\llvert  x-y\rrvert   =1$, this is simply Proposition 32.4 in \cite{spitzer01}. When
$\llvert  x-y\rrvert   >1$, imagine that there are coalescing random walks started at
every position in $\{x, x+1, \ldots, y\}$. Until time $\tau_{x, y}$,
we may regard $X^1$ as the path started from $x$ and $X^2$ as the path
started from~$y$. If all the random walks have coalesced by time $t$,
then $X^1$ and $X^2$ have also coalesced, and hence $\tau_{x,y} \le
t$. Thus,
\[
\prob(\tau_{x,y} > t) \le\sum_{i=x}^{y-1}
\prob( \tau_{i, i+1} >t ) \le\frac{C_0}{\sqrt t} \llvert x - y\rrvert
\]
by the case $\llvert  x-y\rrvert  =1$.\vadjust{\goodbreak}
\end{pf}

\begin{pf*}{Proof of Proposition~\ref{propcomingdownfrominf}}
Suppose there are $m$ distinct particles in the interval $[-Kn, Kn]$ at
time 0. By the pigeonhole principle, there exists at least one pair of
particles that are at a distance of at most $a_{m, n}:= 2Kn/m$. By
Lemma~\ref{lemmahitting}, the probability that these two {\em
unkilled} random walkers will meet each other by time $t_{m, n}:= 4
C_0^2 a_{m, n}^2 $ is at least $1/2$. Hence, in the coalescing system
the probability that there is at most $m-1$ particles after time $
t_{m, n}$ is
certainly at least $1/2$. This happens due to one of the following scenarios:
\begin{longlist}[(iii)]
\item[(i)] At least one of the particles leaves the interval $[ -K n,
Kn]$ and hence gets killed.
\item[(ii)] The two
distinguished particles collide with each other and no others.
\item[(iii)] Some other particle(s) collides with one or both of the
distinguished particles.
\end{longlist}
Moreover, if after time $t_{m, n}$, the number of distinct particles in
the coalescing system still remains $m$ then we can again find at time
$t_{m,n}$ a
possibly different pair of particles that are within distance $a_{m,
n}$ from each other,
and the probability that this pair of particles will collide within the
time interval
$[t_{m,n}, 2t_{m,n}] $ is again at least $1/2$. By repeating this
argument and using the Markov property, we see that if we
let $\tau^{m}_{m-1} = \tau^{m}_{m-1} (n)$ be the first time there are
$m-1$ surviving particles starting from $m$ particles, then, regardless
of the particular initial configuration of the
$m$ particles in $[-Kn, Kn]$,
\[
\prob\bigl(\tau^{m}_{m-1} \ge k t_{m,n} \bigr)
\le2^{-k}.
\]
In particular, $\E[ \tau^{m}_{m-1} ] \le2 t_{m, n}$. Thus, if we
start with one particle at each $x \in[-Kn, Kn] \cap\mathbb Z$, then
the probability that after $\delta n^2$ time the number of particles
remaining is greater than $k$ is, by Markov's inequality, bounded above by
\[
\frac{1}{\delta n^2} \sum_{ m = k+1}^{N} \E
\bigl[\tau ^{m}_{m-1} \bigr] \le\frac{2}{\delta n^2} \sum
_{m=k+1}^{N} t_{m,
n} \le
\frac{32K^2C_0^2}{\delta} \sum_{m=k+1}^{N}
\frac{1}{m^2} \le\frac{32K^2C_0^2}{\delta k},
\]
where $N = \#([-Kn, Kn]\cap\mathbb Z)$. This completes the proof of
the proposition.
\end{pf*}

\subsection{Proof of upper bound}
We now prove (\ref{eUB}).

\subsubsection*{Single tube case} We first prove (\ref{eUB}) in
the case where $m=1$, which is slightly simpler to explain. Set $T =
T_1$, and assume that $\partial_0 T = [a, b] \times\{s\}$. For
$\delta>0$, let $[T_\delta] = [T ] \cap ( \R\times[s+\delta,
t] )$, $\partial_0 T_\delta= [a, b] \times\{s + \delta\}$ and
$\partial_1 T_\delta= \partial_1 T$. Clearly, $T_\delta\in\TUBE
_\square$ for small enough $\delta$ and $T_\delta\le T$. Choose
$K>0$ large enough that $[T] \subseteq[-K/2, K/2] \times\mathbb R$.
For any $a \in\mathbb R$, set $\bar a_\eta= \lfloor a \eta^{-2}
\rfloor\eta^2$.

Since we only care about crossing of the tube $T$, we can work with
coalescing random walks on $\LL_\eta$ that start from every point in
$ I_\eta= [-K,K] \cap\sigma^{-1} \eta\Z$ at time $\bar s_\eta$
and are killed upon leaving $[-K, K]$.
Let $U_\eta$ be the number of distinct particles in the system at time
$\overline{(s+\delta)}_\eta$, and let $z_1^\eta, z_2^\eta, \ldots, z_{U_\eta}^\eta$ denote the space--time positions on $\LL_\eta$ of
these particles at time $\overline{(s+\delta)}_\eta$, enumerated in
some predetermined order. Let $\Gamma^\eta(z_1^\eta, \ldots,
z_{U_\eta}^\eta)$ denote the system of coalescing random walks
started from these space--time positions.
Observe that if $T$ is crossed by the system of coalescing random walks
then necessarily $T_\delta$ is crossed by $\Gamma^\eta(z_1^\eta,
\ldots, z_{U_\eta}^\eta)$.

Now, for all $\eps>0$, by Proposition~\ref
{propcomingdownfrominf}, we can find $k$ so that $\prob( U_\eta>k
) \le\eps$.
Set $\ell= \limsup_{\eta\to0} \P^\eta( \boxminus_T)$, and
assume that this limsup is achieved along a particular subsequence
which we will still denote by $\eta$ with a small abuse of notation.
Then by compactness of $[-K,K] \times[s,t]$, we can find a further
subsequence (also denoted by $\eta$) such that $U_\eta\to U$ and
$(z_1^\eta, z_2^\eta, \ldots, z_{k \wedge U_\eta}^\eta) \to(z_1,
\ldots, z_{k \wedge U})$ in distribution. Along this particular subsequence,
\begin{eqnarray*}
\P^\eta( \boxminus_T) &\le&\prob\bigl( T_\delta
\mbox{ is crossed by } \Gamma^\eta\bigl(z_1^\eta,
\ldots, z_{U_\eta}^\eta\bigr)\bigr)
\\
& \le& \eps+ \prob\bigl(T_\delta\mbox{ is crossed by }
\Gamma^\eta \bigl(z_1^\eta, \ldots,
z_{k \wedge U_\eta}^\eta\bigr) \bigr)
\\
& \to& \eps+ \prob\bigl( T_\delta\mbox{ is crossed by }
\Gamma(z_1, \ldots, z_{k \wedge U})\bigr)
\end{eqnarray*}
by Lemma~\ref{lkpoint}, the Markov property and the bounded
convergence theorem where $ \Gamma(z_1, \ldots, z_{k \wedge U})$
denotes coalescing Brownian motions started from $(z_1, \ldots, z_{k
\wedge U})$. By Theorem~\ref{tBW}, we conclude
\[
\ell\le\eps+ \P_\infty( \boxminus_{T_\delta}). %
\]
Since $\eps$ is arbitrary, $\ell\le\P_\infty( \boxminus
_{T_\delta})$. Now, as $\delta\to0$, the events $\boxminus
_{T_\delta}$ are decreasing, so
\[
\lim_{\delta\to0} \P_\infty( \boxminus_{T_\delta}) =
\P_\infty \biggl( \bigcap_{\delta>0}
\boxminus_{T_\delta} \biggr).
\]
Since $\HH$ consists of \emph{closed} collection of tubes, $ \bigcap_{\delta>0} \boxminus_{T_\delta} = \boxminus_T$. Thus,
%
\begin{equation}
\label{econclm=1} \ell\le\P_\infty( \boxminus_{T})
\end{equation}
and so (\ref{eUB}) holds in the case $m=1$.

\subsubsection*{Multi tube case}
Now assume that $m \ge1$, and to keep notation simple we will assume
that $m=2$. Let $T, T'$ be two tubes in $\TUBE_\square$ with the
start and the end times $s, t$ and $s', t'$, respectively. We assume
without loss of generality that $s< s'$. Reasoning as in the case $m
=1$, it is easy to deal with the case where $[T] \cap[T'] = \varnothing
$. We thus assume that $[T] \cap[T'] \neq\varnothing$, and hence $s' < t$.
For $\delta>0$ small enough, let $ T_\delta$ and $T'_\delta$ be
tubes in $ \TUBE_\square$ obtained similarly from $T$ and $T'$ as in
the single tube case.
Further, the tube $ T_\delta$ is decomposed into two tubes
$T^-_{\delta}$ and $T^+_{\delta}$ in $\hat\TUBE$ such that
\[
\bigl[T^-_{\delta} \bigr] = [T ] \cap \bigl(\R\times\bigl[s + \delta,
s'+ \delta \bigr] \bigr) \quad\mbox{and}\quad \bigl[T^+_{\delta}
\bigr] = [T ] \cap \bigl( \R\times\bigl[s' + \delta, t\bigr] \bigr).
\]

Choose $K>0$ large enough that $[T], [T'] \subseteq[-K/2, K/2] \times
\R$. Essentially, we wish to consider coalescing random walks in $\LL
_\eta$ that start from every point in $ I_\eta= [-K,K] \cap\eta\Z$
at time $\bar s_\eta$ and $\bar{s}'_\eta$, that are killed upon
leaving $[-K, K]$. It is useful to picture the particles starting at
time $\bar s_\eta$ as being colored \emph{blue} and those starting at
time $\bar s'_\eta$ as being colored \emph{red}. We wish to apply the
reasoning of the case $m=1$ separately to all three tubes above, but we
need to be a little careful to avoid interactions between the blue and
red particles during the interval $[s', s'+\delta]$.

In order to do so, we introduce a coupling of red and blue particles
which dominates the coalescing random walks. Here is the precise
definition. We start by associating to each vertex $z\in\LL_\eta$ a
random variable $\xi_z$ which is an i.i.d. copy of the step
distribution $\xi$. Coalescing random walk $(S_k, t_k)_{k \ge0}$ in
$\LL_\eta$ from space--time point $z = (x_0, t_0) \in\LL_\eta$ may
be defined by setting, $S_0 = x_0$ and for $k \ge0$,
%
\begin{equation}
\label{ecrw} S_{k+1} = S_{k} + \sigma^{-1} \eta
\xi_{(S_{k-1}, t_{k-1})}\quad \mbox{and} \quad t_k = t_0 + k
\eta^2.
\end{equation}

In order to prevent the blue and red particles from interacting during
the interval $[\overline{s'}_\eta, \overline{(s'+\delta)}_\eta]$,
we modify this description as follows. Consider all the points $z \in
I_\eta\times[\overline{s'}_\eta, \overline{(s'+\delta)}_\eta]
\cap\LL_\eta$.
We endow each such $z$ with a new independent copy $\xi'_z$ of $\xi$
in addition to the original $\xi_z$. Then a blue particle will use the
random variable $\xi_z$ to move forward from space--time point $z$, but
a red particle will always use the random variable $\xi'_z$ if it has
the choice between $\xi_z$ and $\xi'_z$ (otherwise it uses $\xi_z$).
The particles still get killed upon exiting the interval $[-K, K]$.
Note that as a result of this definition, after time $\overline
{(s'+\delta)}_\eta$, if a red and a blue particle are on the same
site then they coalesce and necessarily follow the same path afterward.
We may thus think of the resulting particle as carrying both the red
and blue colours (see Figure~\ref{fcoupling}).

%
\begin{figure}

\includegraphics{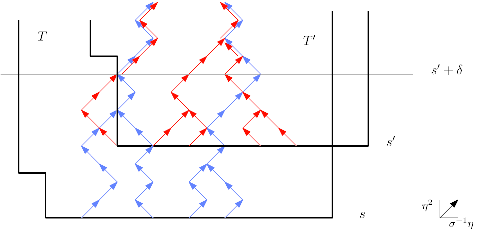}

\caption{The coupling of partially coalescing (simple) random walks.}
\label{fcoupling}
\end{figure}

Then observe that if $\boxminus_{T} \cap\boxminus_{T'}$ occurs
then necessarily the following three things must occur:
\begin{longlist}[(a)]
\item[(a)] $T_{\delta}^-$ is crossed by the blue particles.
\item[(b)] $T_{\delta}^+$ is crossed by the remaining blue particles
at time $\overline{(s'+\delta)}_\eta$. 
\item[(c)] $T'_\delta$ is crossed by the remaining (blue or red)
particles at time $\overline{(s'+\delta)}_\eta$.
\end{longlist}

Let $U_\eta^{-} $ be the number of distinct blue particles in the
system at time $\overline{(s+\delta)}_\eta$, and let $z_1^\eta,
z_2^\eta, \ldots, z_{U_\eta^-}^\eta$ denote the space--time positions
on $\LL_\eta$ of these particles. Let $\Gamma_-^\eta= \Gamma
_-^\eta(z_1^\eta, \ldots, z_{U_\eta^-}^\eta)$ denote the system of
coalescing random walks started from these space--time positions, and
ended at time $\overline{(s'+\delta)}_\eta$. Let $U_\eta^+$ denote
the number of blue particles in $\Gamma_-^\eta(z_1^\eta, \ldots,
z_{U_\eta-}^\eta)$ left at time $\overline{(s'+\delta)}_\eta$.
Let also $U_\eta'$ denote the number of distinct red particles at time
$\overline{(s'+\delta)}_\eta$. Denote by $w_1^\eta, \ldots, w^\eta
_{U^+_\eta}, y_1^\eta, \ldots,  y_{ U'_\eta}^\eta$ the locations of
the blue and red particles at time $\overline{(s'+\delta)}_\eta$,
and denote by $\Gamma_+^\eta= \Gamma_+^\eta(w_1^\eta, \ldots,
w^\eta_{U^+_\eta}, y_1^\eta,\break \ldots, y_{ U'_\eta}^\eta)$ the
collection of coalescing blue and red particles started from these
positions at time $\overline{(s'+\delta)}_\eta$. Note that blue and
red particles are allowed to coalesce after this time (in which case
they count as both blue and red particles for what follows).

From the above discussion, it follows that
%
\begin{eqnarray}\label{econd}
\P^\eta( \boxminus_{T} \cap
\boxminus_{T'})
&\le&\prob\bigl(T_\delta ^-, T_\delta^+
\mbox{ are crossed by } \Gamma_-^\eta\mbox{ and }
T'_\delta\mbox{ is crossed by } \Gamma_+^\eta
\bigr).\hspace*{-25pt}
\end{eqnarray}
\margin{C: I changed the above equation according to the new (a),(b),(c)}
As in the case $m=1$ set $\ell= \limsup_{\eta\to0} \P^\eta(
\boxminus_{T} \cap\boxminus_{T'})$, and assume that this limsup is
achieved along a particular subsequence which we will still denote
$\eta$ with a small abuse of notation. For all $\eps>0$, by
Proposition~\ref{propcomingdownfrominf}, we can find $k$ so that
$\prob(\max( U_\eta^-, U_\eta') >k ) \le\eps$ for all $\eta$.
Then by compactness, we can find a further subsequence (also denoted by
$\eta$) such that $ ( U_\eta^-, U_\eta') \to( U^-, U')$ and the two
vectors $(z_1^\eta, \ldots, z_{k \wedge U_\eta^-}^\eta)$ and
$(y_1^\eta, \ldots, y_{U'_\eta\wedge k}^\eta)$ converge jointly to
two vectors
$(z_1, \ldots, z_{k \wedge U^-})$ and $(y_1, \ldots, y_{k \wedge
U'})$ in distribution. Let $\Gamma(z, y)$ denote coalescing Brownian
motions starting from $(z_1, \ldots, z_{k \wedge U^-},y_1, \ldots,
y_{k \wedge U^+})$.
Applying the Markov property in (\ref{econd}) repetitively and by
bounded convergence theorem, we see after taking the limit that along
this particular subsequence:
\begin{eqnarray*}
\ell& \le&\eps+ \prob\bigl( T_\delta, T'_\delta
\mbox{ are crossed by } \Gamma(z,y)\bigr)
\\
& \le&\eps+ \P_\infty(\boxminus_{T_\delta} \cap\boxminus
_{T'_\delta}).
\end{eqnarray*}
From there, we conclude as in (\ref{econclm=1}), since $\eps>0$ is
arbitrary and $\HH$ consists of closed collections of tubes,
\[
\ell\le\P_\infty(\boxminus_T \cap\boxminus_{T'}).
\]
This proves (\ref{eUB}), and hence Theorem~\ref{tscalinglimit}.
\end{pf*}

\section{Coalescing flow on Sierpinski gasket}\label{Sgasket}
\subsection{Sierpinski gasket}
Let $H_0$ be the unit triangle in $\mathbb R^2$ with vertices $ \{
(0,0), (1, 0), (1/2, \sqrt3 /2) \}$.
The {\em finite Sierpinski gasket} is
a fractal subset of the plane that can be constructed via the following
Cantor-like cut-out procedure. Let $\{ b_0, b_1, b_2\}$ be the
midpoints of three sides of $H_0$ and let $A$ be the interior of the
triangle with vertices $\{b_0, b_1, b_2\}$. Define $H_1:= H_0
\setminus A$ so that $H_1$ is the union of $3$ closed upward facing
triangles of side length $2^{-1}$. Now repeat this operation on each of
the smaller triangles to obtain a set $H_2$, consisting of $9$ upward
facing closed triangles, each of side $2^{-2}$. Repeating this
procedure, we have a decreasing sequence of closed nonempty sets $\{
H_n\}_{n=0}^\infty$ and we define the finite Sierpinski gasket as
\[
G_{\mathrm{fin}}:= \bigcap_{n=0}^\infty
H_n.
\]
We call the unbounded set
\[
G:= \bigcup_{n=0}^\infty2^n
G_{\mathrm{fin}}
\]
the {\em infinite Sierpinski gasket}. We equip it with the shortest
path metric $\rho_G$ which is comparable to the usual Euclidean metric
$\llvert   \cdot\rrvert  $ (see, e.g., \cite{barlow98}, Lemma 2.12) with the relation,
\[
\llvert x- y\rrvert \le\rho_G(x,y) \le c\llvert x-y\rrvert \qquad
\forall x, y \in G,
\]
for a suitable constant $1 < c < \infty$. Let $\mu$ denote the
$d_f$-dimensional Hausdorff measure
on $G$ where $d_f:= \log3/ \log2$
is the {\em fractal} or {\em mass dimension} of the gasket.
The following estimate on the volume growth of $\mu$ is known (see
\cite{barlow88}):
%
\begin{equation}
\label{eqvolgrowth} \mu\bigl( B(x, r) \bigr) \le C r^{d_f}\qquad\mbox{for }
x \in G, 0 < r < 1,
\end{equation}
where $B(x,r) \subseteq G$ is the open ball with center $x$ and radius $r$
in the Euclidean metric.

For each $n \in\Z$, the set $2^m H_{m+n}$ is made up of $3^{m+n}$
triangles of side length $2^{-n}$ whenever $m \ge n$. Each of those
triangles are called an {\em$n$-triangle} of $G$.
Denote by $\cS_n$ the collection of all $n$-triangles of $G$.
Let $\cV_n$ be the set of vertices of the $n$-triangles. We will
restrict our attention to infinite Sierpinki gasket while constructing
the coalescing Brownian flow. The case of finite gasket can also be
dealt with very similar arguments.

\subsection{Brownian motion on gasket}

We construct a graph $G_n$ embedded in the plane
with vertices $\cV_n$
by adding edges between pairs of vertices that are
distance $2^{-n}$ apart from each other. Let $X^n$ be the
nearest-neighbor random walk on $G_n$ simultaneously defined on the
same probability space. It is known (see \cite{barlow88,barlow98}) that
the sequence $(X^n_{\lfloor5^nt \rfloor})_{ t \ge0}$ converges
almost surely as $n \rightarrow\infty$
to a limiting process $(X_t)_{t \ge0}$ that is a $G$-valued strong
Markov process (indeed, a Feller process)
with continuous sample paths. The process $X$ is naturally called the
Brownian motion on the gasket. It has the following scaling property:
%
\begin{equation}
\label{ebmscaling} (2 X_t)_{t \ge0} \mbox{ under }
\prob^x\mbox{ has same law as }(X_{5t})_{t \ge0}
\mbox{ under } \prob^{2x}.
\end{equation}
The process $X$ has a symmetric transition density $p_t(x, y), x, y \in
G, t>0$ with respect to the measure $\mu$ that is jointly continuous
on $(0, \infty) \times G \times G$.
Let $ d_w:= \log5 / \log2$ denote the {\em walk dimension}
of the gasket. The following crucial ``heat kernel bound'' is
established in \cite{barlow88}
%
\begin{equation}
\label{eqdensitybound} p_t( x, y) \le c_1 t^{ -d_f/ d_w} \exp
\biggl( - c_2 \biggl( \frac
{\llvert  x-y\rrvert  ^{d_w}}{t} \biggr)^{ 1/ ( d_w-1)}
\biggr).
\end{equation}
A matching lower bound (with different constants $c_3$ and $c_4$) also
exists. This shows that the Brownian motion on the gasket is sub-diffusive.

\subsection{Coalescing Brownian flow on gasket}


We now state the analogue of Theorem~\ref{tBW} in the case of the
Sierpinski gasket.
Let $\cD=\{z_1, z_2, \ldots\}$ be a countable ordered set which is
dense in $G\times\R$ where $z_i = (x_i, t_i)$. Let $(W_j)_{j \ge1}$
be an independent family of Brownian motions on the Sierpinski gasket
$G$ defined on a common probability space $(\Omega, \mathcal{F},
\prob)$, started from time $t_j$ at position $x_j$. We can apply the
coalescing rule to obtain a collection $\Gamma_n$
of $n$ coalescing Brownian motions on $G$ denoted by $W^c_1, W^c_2,
\ldots, W^c_n$. Let us define $\calW_n$ to be, as before, the set of
all tubes (now $d=2$) crossed by $\Gamma_n $, that is, $\calW_n =
\calW(z_1, \ldots, z_n): = \Cr(\Gamma_n) \in\HH$.
The next theorem defines the coalescing Brownian flow on the gasket and
its proof is a straightforward adaptation of the arguments in the proof
of Theorem~\ref{tBW} (we leave the details to the reader).

\begin{theorem}\label{tBWg}
As $n \to\infty$, $\cW_n$ converges in distribution to a random
variable $\cW_\infty$, whose law does not depend on $\cD$ (including
its order).
\end{theorem}

\begin{definition}
\label{dBWg} A random variable on $\HH$ with law $\P_\infty$ is
called a {\em coalescing Brownian flow} on the Sierpinski gasket.
\end{definition}


\subsection{Characterization}

In a way which is analogous to Theorem~\ref{tBW2}, we state a useful
characterization of the coalescing Brownian flow on the Sierpinski
gasket. The result will be formally very similar to Theorem~\ref
{tBW2} but we will work with a slightly different class of tubes,
chosen so that they are more suited to the\vspace*{1pt} geometry of the Sierpinski
gasket. Let $\tr_0$ denote the convex hull of the vertices $0, 1$ and
$e^{i \pi/ 3}$ in $\R^2$. Let $\mathbb{T}_n$ be the triangular
lattice on the plane with mesh size $2^{-n}$ (so that $G_n$ is a
subgraph of $\mathbb{T}_n$). Define
\[
\cE= \biggl\{z+ \frac{1}{2^{n}} \tr_0\dvtx  n \in\Z, z \in
\mathbb{T}_n \biggr\}.
\]
Note that if $\tr\in\cE$, $\tr$ is an upward-facing equilateral
triangle in $\mathbb{T}_n$ for large enough $n$. The Brownian motion
starting at some point inside $\tr$ can escape $\tr$ only through one
of the three vertices of $\tr$.

\begin{definition}[(A dense family of triangular tubes)]\label{dktubesg}
Let $\triT$ be the family of all tubes $T$ such that:
\begin{longlist}[(a)]
\item[(a)] The set $[T]$ can be expressed as an union of a {\em
finite} number of cylinders (triangular prisms) of the form $\tr\times
[s,t]$ with $s<t$, for some $\tr\in\cE$.
\item[(b)] $\partial_0 T = [T] \cap(\R^2 \times\{t_0\})$ and
$\partial_1 T = [T] \cap(\R^2 \times\{t_1\})$, where $t_0$ and
$t_1$ are the start and the end time of $T$.
\end{longlist}
\end{definition}

One can check that $\triT$ is dense in $\TUBE$ because every
downward-facing triangle is a subset of an upward-facing triangle of
twice its size.


We now state our characterization.

\begin{theorem}\label{tBW2g}
$\P_\infty$ is the unique probability measure on $\HH$ such that for
all $m \ge1$ and for all fixed tubes $T_1, \ldots, T_m \in\triT$,
%
\begin{eqnarray}\label{eqBWg}
&& \P_\infty(\boxminus_{T_1}\cap\cdots\cap
\boxminus_{T_m} )
\nonumber\\[-8pt]\\[-8pt]\nonumber
&&\qquad  = \sup_{{n \ge1; z_1, \ldots, z_n \in G\times\R}} \prob \bigl( \calW
(z_1, \ldots, z_n) \in\boxminus_{T_1}\cap\cdots
\cap\boxminus _{T_m} \bigr).
\end{eqnarray}
\end{theorem}

\begin{pf}
The proof of Theorem~\ref{tBW2} can be adapted as follows. Lemma~\ref
{Lcdi}, which relies on the fact that Arratia's flow comes down from
infinity, will now use Theorem 5.1 in \cite{ems09} which says that the
same is true even for coalescing Brownian motion on the gasket.
But Lemma~\ref{LA2} needs a different argument, since it relies on
scale invariance of (real) Brownian motion at all scales. Since
Brownian motion on the gasket is scale-invariant only for a discrete
set of scales (which does not come arbitrarily close to 1), this lemma
needs a different proof.

Thus, let $T \in\triT$ be fixed and let $\partial_0 T = B \times\{
t_0\}$, where $B$ is the union of finitely many triangles from $\cE$.
To keep the presentation simpler we will assume that $B = \tr$ for
some fixed equilateral triangle $\tr$ in $\cE$ of side $1$. Let
$t_e>0$ be such that $T \cap(\R^2 \times[ t_0, t_0 + t_e]) = \tr
\times[t_0, t_0 + t_e]$. Let $t_1 = t_0 + \delta$ and $t_2 = t_1+
\alpha$, where $\alpha= t_e/4$. Fix a countable dense set $\cD$ of
$G \times\R$ and consider a countable system of coalescing Brownian
particles starting from the space--time points in $\cD$.
Let $A_\tr(s,t)$ denote the set of locations at time $t$ of the
particles that started at some time before $s$ and were in $\tr$
throughout the time interval $[s,t]$. For $\delta>0$, let $\tr^\delta
$ be the $\delta$-enlargement of $\tr$. We wish to show that for all
$\eps>0$, we can choose $\delta>0$ so that $A_\tr(t_1, t_2) = A_{\tr
^{\delta}}(t_0, t_2)$ with probability greater than $1- \eps$. Again,
by continuity of the distribution of $\llvert  A_\tr(s, t)\rrvert  $ in $s $ and $t$,
it suffices to show the following.

\begin{lemma} \label{Lenlarg} Given $\eps>0$, there exists $\delta
_0 >0$ such that for all $\delta\le\delta_0$, the event
$\llvert  A_\tr(t_1, t_2)\rrvert   = \llvert  A_{\tr^{\delta}}(t_1, t_2)\rrvert  $ holds with
probability greater than $1-\eps$.
\end{lemma}

\begin{pf}
Let $v_1, v_2, v_3$ be the vertices of the triangle $\tr$. Call the
vertex $v_i $ an \textit{exit point} \margin{Added boldface here} of $\tr
$ if the set $B(v_i, \lambda) \setminus\tr$ has nonempty
intersection with $G$ for each $\lambda>0$. Similarly, call the vertex
$v_i $ an \textit{entry point} of $\tr$ if the set $B(v_i, \lambda)
\setminus\{ v_i\}$ has nonempty intersection with $\tr\cap G$ for
each $\lambda>0$.
For the vertex $v_i$, let $q_k^i $ (resp., $r_k^i$) be the union of
one or two (resp., one) $k$-triangle(s) in $\cS_k$ attached to $v_i$
and lying outside (resp., inside) of $\tr$, if $v_i$ is an exit
(resp., entry) point of $\tr$ and the empty set otherwise (see
Figure~\ref{fgasket}).
Since $\tr$ is an upward-facing triangle adapted to the triangular
lattice $\mathbb{T}_0$, the difference between $\tr^\delta\cap G$
and $\tr\cap G$ is ``tiny'' when $\delta$ is small. Indeed, for $\delta
= 2^{-(k+1)}$, $(\tr^\delta\setminus\tr) \cap G \subseteq q_k^1
\cup q_k^2 \cup q_k^3$. Thus, we can assume that each $v_i$ is an exit
point of $\tr$, if not we can safely ignore it for the rest of the
proof. For $\delta= 2^{-(k+1)}$, let $\tr^{-\delta}$ be the closed
set obtained by removing $r^1_k, r^2_k$ and $r^3_k$ from $\tr$. We
first show that:

\begin{claim}\label{claim1}
One can choose $\delta_1>0$ small enough
that with probability at least $1-\eps/2$, no particle starting in
$\tr^\delta\setminus\tr^{-\delta}$ at time $t_1$ can stay inside
$\tr^\delta$ up to time~$t_2$, for any $\delta< \delta_1$.
\end{claim}

The proof of the above claim relies on the following further claim.

%
\begin{figure}

\includegraphics{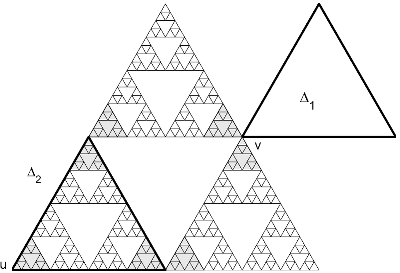}

\caption{The triangles $\tr_1, \tr_2 \in\cE$ are shown with thick lines.
For example, $u$ is an entry point of $\tr_1$ but not an exit point
where as $v$ is an exit point of $\tr_2$ but not an entry point. The
shaded region attached to $v$ is $q_k^v$ and the shaded triangle
attached to $u$ is $r_k^u$ (for $k=3$).}
\label{fgasket}
\end{figure}

\begin{claim}\label{claim2}
Given $\eps_1 >0$ and $\gamma>0$, there
exists $k_0$, such that for all $k \ge k_0$, the total coalescence time
of the countable particles starting from {\em any} $k$-triangle of $G$
is less than $\gamma$ with probability $1 - \eps_1$.
\end{claim}

\begin{pf*}{Proof of Claim \ref{claim1}} To see this, note that for $\delta=
2^{-(k+1)}$, the set $\tr^\delta\setminus\tr^{-\delta}$ consists
of three connected regions, say $\ell^1_k, \ell^2_k$ and $\ell^3_k$,
of diameter bounded by $2^{-(k-1)}$ where $\ell^i_k = q^i_k \cup
r_i^k$. Consider the particles starting from $\ell_k^i$, $1\le i \le
3$ at time $t_1$. By Claim~\ref{claim2}, for $\gamma>0$ fixed but sufficiently
small, at time $t_1 + \gamma$, the number of descendants of these
particles is one with probability at least $1 - \eps/4$, and hence
coincides with the particle starting from $v_i$ at time $t_1$. By
choosing $\gamma>0$ small and $\delta>0$ even smaller, the particle
starting from $v_i$ at time $t_1$ is guaranteed to exit the triangle
$q_k^i$ via one of its vertices other than $v_i$ between time $t_1+
\gamma$ and $t_2$ with probability $1-\eps/4$. On the intersection of
these two events (which has probability at least $1- \eps/2$), no
particle starting from $\ell^i_k$ can cross $\tr^\delta\times[t_1, t_2]$.
\end{pf*}

\begin{pf*}{Proof of Claim~\ref{claim2}}
We first claim that there exist $n \ge1$ and $u\in(0,1)$, depending
only on $\eps_1$, such that with probability at least $1 - \eps_1/2$,
the number of descendants at time $u$ of the particles starting from
any $0$-triangle $L$ in $\cS_0$ is less than $n$ and moreover, those
particles are inside $L+ B(0, 1/2)$ at time $u$. \margin{C:
corrected as suggested by Arnab. Please check, A: Checked.} The main
content of the claim lies in the fact that $n$ and $u$
can be chosen independent of $L$. This follows from Lemma~5.6 of \cite
{ems09}, which actually gives much more in terms of quantitative
bounds. Thus, by scaling, for any $k \ge1$, with probability at least
$1 - \eps_1/2$, the number of descendants at time $u  5^{-k}$ of
particles starting from any \mbox{$k$-}triangle $L$ is less than $n$, and
moreover, those particles are inside $L+ B(0, 2^{-(k+1)})$ at time $u
5^{-k}$. Now choose $k_0$ large enough such that for all $k \ge k_0$
two independent Brownian motions starting within distance $2^{-(k-1)}$
from one another at time $0$ will meet by time $\gamma/2$ with
probability at least $1 - \eps_1/(2n)$. This is possible by (\ref
{B2pt}) whose proof is given in Section~\ref{subsecinvgasket}. This
completes the proof of Claim~\ref{claim2}.
\end{pf*}


Now we will continue with the proof of Lemma~\ref{Lenlarg}. Consider
the particles starting from $\tr^{-\delta_1}$ at time $t_1$. We have
to argue that for sufficiently small $\delta$ with probability $1-
\eps/2$, no such particle can stay within $\tr^\delta$ between time
$[t_1, t_2]$ but go outside $\tr$ at some time in $[t_1, t_2]$. Let
$\Xi_t$ be the closure of the locations of the coalescing particles at
time $ t_1 + t$ which were in the compact set $\tr^{-\delta_1}$ at
time $t_1$. It was proved in \cite{ems09}, Theorem 5.2, that $\Xi_t
\to\Xi_0 = \tr^{-\delta_1}$ in probability as $t \to0+$ in the
Hausdorff metric.
Consequently, we can find $\eta>0$ sufficiently small such that no
particle which was inside $\tr^{-\delta_1} $ at time $t_1$ can leave
$\tr$ before time $t_1+ \eta$, with probability greater than $1-\eps
/6$. By coming down from infinity, choose $n$ sufficiently large, so
that $\llvert  A_{\tr}(t_1, t_1+\eta) \rrvert  \le n$ with probability at least
$1-\eps/6$.

Now the lemma will follow if we show that for any particle, say $X$,
which is inside $\tr$ at time $t_1 +\eta$, the probability that it
leaves $\tr$ some time between $t_1+ \eta$ and $t_2$ but always stays
inside $\tr^{\delta}$ within time $[t_1+\eta, t_2]$ is at most $\eps
_0 = \eps/(6n)$. Let $\tau$ be time when $X$ hits one of the vertices
of $\tr$ for the first time after $t_1 +\eta$. We find $\beta\in(0,
t_2 - t_1 - \eta)$ small such that $\prob( \tau\in[t_2-\beta,t_2
]) \le\eps_0/2 $, uniformly over the position of $X$ in $\tr$ at
time $t_1+\eta$. If $\tau\le t_2 - \beta$, then $X_\tau$ will be at
one of the vertices $v_i$ at time $\tau$ and by choosing $\delta$
small, it will leave $\tr^\delta$ some time during the time interval
$[\tau, \tau+ \beta]$ with probability at least $1 - \eps_0/2$.

Putting together these observations we have completed the proof of the lemma.
\end{pf}

The rest of the proof of Theorem~\ref{tBW2g} is a straightforward
adaptation of Lemmas~\ref{LA1} and~\ref{lemcorner}, which is left
to the reader. Note in particular that for a tube $T \in\triT$ with
$[T] = \bigcup_{j=1}^k \tr_j \times[s_j, t_j]$, the set\vspace*{1pt} of corner
points of $T$ is a subset of $\{ (v_{ij}, s_j), (v_{ij}, t_j)\dvtx  1 \le i
\le3, 1 \le j \le k \}$, and hence is finite. Here, $v_{1j}, v_{2j}$
and $v_{3j}$ denote the vertices of the triangle $\tr_j \in\cE$.
\end{pf}




\section{Invariance principle for coalescing random walks on gasket}\label{Sgasket2}

Let $\eta= 2^{-n}$ and consider an infinite Sierpinsky gasket $G_n$
with mesh size $2^{-n}$. In this section, we consider coalescing random
walks on $G_n$ defined as follows: initially there is a particle at
every vertex of $G_n$. They perform independent simple random walks,
jumping every $5^{-n}$ units of time, and coalesce when they are on the
same vertex of $G_n$. Consider the law $\P^\eta$ on $\HH$ that these
particles induce.
\margin{C: similarly as for RWs on $\Z$, one needs to address the
aperiodicity issue of the walks. Fortunately, it is the case, but where
does this implicitly appear in the proof ? Ok, found it. I added a
remark below.}

\begin{theorem}\label{tscalinglimitg}
As $\eta\to0$,
\[
\P^\eta{\to} \P_\infty,
\]
weakly, where $\P_\infty$ is the law of the coalescing Brownian flow
on $G$, as defined in Theorem~\ref{tBWg}.
\end{theorem}

\begin{pf}
The proof of Theorem~\ref{tscalinglimitg} follows the same outline
as in the Brownian case. However, two important facts are needed to
make the strategy applicable. The first is a statement about ``uniform
coming down from infinity'' (stated below in Proposition~\ref
{PCDIs}). The second is an intuitively obvious statement that
finitely many coalescing random walks converge to the same number of
coalescing Brownian motions. This will be stated in Proposition~\ref
{pkpoint}. This is probably well known in the folklore, but we could
not find a reference for it.

\subsection{Uniform coming down from infinity on the gasket}

Consider the setup above, with coalescing random walks started at time
0 on $G_n$, and jumping every $5^{-n}$ units of time. Give a bounded
region $T \subset\R^2$, suppose that the random walk particles are
killed as soon as they touch $\R^2 \setminus T$. Let $N(t)$ denotes
the number of particles left at time $t$. The following is the analogue
of Proposition~\ref{propcomingdownfrominf} but for the case of the
gasket. As before, this is essentially the only place where one needs
some quantitative estimates about coalescing random walks.

\begin{proposition}[(Uniform coming down from infinity)]\label{PCDIs}
For every $\delta>0$ and $\eps>0$, there exists $k \ge0$ depending
only on $T$, such that $\prob(N(\delta) > k) \le\eps$ for all
sufficiently large $n$.
\end{proposition}

\begin{pf}
The proof depends on two simple lemmas. In the rest of the proof, we
assume without loss of generality that $T$ is the unit equilateral
triangle with apices at $z=0$, $z=1$, and $z = e^{i\pi/ 3}$ when
viewed as a subset of the complex plane.
\end{pf}

\begin{lemma}\label{Lcoaltimepair}
There exists a universal constant $0<C<\infty$ such that the following
holds for all $1\le k \le n$. Let $X,Y$ be two independent (unkilled)
random walks started at $x,y \in G_n$ such that $\llvert  x-y\rrvert   \le2^{-k}$. Then
\[
\prob\bigl( \tau< 5^{-k}\bigr) \ge\frac{1}{C} %
\]
uniformly on $x,y\in G_n$ such that $\llvert  x-y\rrvert   \le2^{-k}$, where $\tau$
is the first meeting time of $X$ and $Y$.
\end{lemma}

\begin{pf}
This relies on uniform heat-kernel estimates of the random walks on the
gasket due to Jones \cite{jones96}.
\margin{I guess Jones handles the possible issue of aperiodicity here,
I will add a remark, but please check.}
To use this, it is convenient to use a different scaling of space and
time: thus consider $G_0 = 2^n G_n$ and let the random walks $X,Y$ make
jumps at integer times. Then Lemma~\ref{Lcoaltimepair} is equivalent
to the statement that, uniformly over $x,y \in G_0$ with $\llvert  x-y\rrvert   \le
2^k$, $\prob(\tau< 5^k) \ge1/C$. Let $p_t(x, y) = \prob(X_t = y\mid X_0
= x)$ be the transition density function of the random walk $X$.
Now, Theorem 18 of \cite{jones96} states that for all $t \ge c_0 \llvert  x-y\rrvert
\vee t_0$,
%
\begin{equation}
\label{heatlb} p_t(x,y)\ge c_1 t^{-d_s/2} \exp
\biggl(- c_2 \biggl(\frac
{\llvert  x-y\rrvert  ^{d_w}}{t} \biggr)^{\afrac{1}{d_w -1}} \biggr),
\end{equation}
where
%
\begin{equation}
\label{dimensions} d_s = \frac{2\log3}{\log5};\qquad d_w =
\frac{\log5}{\log2}.
\end{equation}
Likewise, Theorem 17 of \cite{jones96} states that for all $t \ge c_0
\llvert  x-y\rrvert   \vee t_0$, \margin{C: added $c_0$ here}
%
\begin{equation}
\label{heatub} p_t(x,y)\le c_3 t^{-d_s/2} \exp
\biggl(- c_4 \biggl(\frac
{\llvert  x-y\rrvert  ^{d_w}}{t} \biggr)^{\afrac{1}{d_w -1}}
\biggr).
\end{equation}
Let $J$ denote the number of intersections of the walks $X,Y$ during
the time interval $[1, 5^k]$. Then by reversibility, and (\ref
{heatlb}), for $k$ large, \margin{added $=z$ below}
\begin{eqnarray*}
\E[J] & =& \sum_{s=1}^{5^k} \sum
_{z \in G_0} \prob(X_s = Y_s=z) = \sum
_{s=1}^{5^k} p_{2s}(x,y)
\\
& \ge& c_1 \sum_{s= (1/2) 5^k}^{5^k}
\frac{1}{(5^k)^{\log3/ \log5}} \exp \biggl( - c_2 \biggl(\frac{(2^k)^{\log5 / \log2}}{5^k}
\biggr)^{1/(d_w -1)} \biggr) \ge c (5/3)^k.
\end{eqnarray*}
On the other hand, by (\ref{heatub}), for $k$ large,
\begin{eqnarray*}
\E\bigl[J^2\bigr] 
& \le&\E[J] + 2 \sum
_{1\le s < t \le5^k} \prob(X_s = Y_s) \sup
_{w
\in G_0} p_{2(t-s)}(w,w)
\\
& \le&\E[J] + 2 \sum_{s=1}^{5^k}
\prob(X_s = Y_s) \sum_{s< t \le
5^k}
\frac{c'}{(t-s)^{ \log3 / \log5}}
\\
& \le&\E[J] + c' \sum
_{s=1}^{5^k} \prob(X_s =
Y_s) (5/3)^k \le c' (5/3)^k
\E[J].
\end{eqnarray*}
Thus, by the Payley--Zygmund inequality,
\[
\prob(J>0) \ge\frac{\E[J]^2}{ \E[J^2]} \ge\frac{\E[J] \cdot
c(5/3)^k}{c' (5/3)^k \E[J] } \ge c/c',
\]
as required.
\end{pf}

\begin{remark}\label{raperiodic}
Note that Lemma~\ref{Lcoaltimepair} implicitly relies on the fact
that the simple random walk on $G_n$ is aperiodic [similarly as the
distribution $\xi$ in~(\ref{axi})]. Though we do not directly appeal
to aperiodicity in our proof, it is used in Jones' lower bound estimate~(\ref{heatlb}). \margin{C: added this remark}
\end{remark}

Our second lemma is an induction scheme which is inspired by an
argument in \cite{ems09} for the fact that coalescing Brownian
particles come down from infinity.
See also \cite{global} where a similar argument is used.

\begin{lemma}
\label{Linduction} Fix $t>0$ and let $m$ be such that $3^{m+1} \le N=
N(t) < 3^{m+2}$. Then there exists an absolute constant $\theta<1$,
such that by time $t' = t + 5^{-m}$, $N(t') \le\theta N(t)$ with
probability greater than $1- e^{- c N}$.
\end{lemma}

\begin{pf}
We first claim that it is possible to find a pairing of the particles
$(X_i, Y_i)_{1 \le i \le N/2}$ such that $\llvert  X_i(t) - Y_i(t)\rrvert   \le2^{-m}$
for all $1\le i \le\tfrac{N}{3}$. Indeed, tile $T$ with $3^{m}$
triangles of side length $2^{-m}$. Within each such triangle, pair as
many particles as possible. This leaves at most one unpaired particle
per triangle, and we pair these arbitrarily.

Having constructed this pairing, consider now a partial coalescing
system in which coalescence occurs only between matched particles, and
distinct pairs of particles evolve completely independently of one
another. By the monotonicity property of coalescing random walks, it is
easy to argue (see \cite{ems09}) that the partial system dominates in
distribution the fully coalescing system. Hence, it suffices to prove
the claim on the partially coalescing system.

By Lemma~\ref{Lcoaltimepair} and by the strong Markov property, at
time $t'$, each pair $(X_i, Y_i)$ with $1\le i \le\tfrac{N}{3}$ has a
probability at least $1/C$ to have coalesced, and these events are
independent of one another. (Note further that the number of particles
may decrease further due to particles leaving the region $T$, but this
can only help us.) Thus, $N(t') \le N - (1/C') (N/3)$ with probability
greater than $1- \exp( - cN')$, by an easy large deviation bound on
the binomial random variables. Taking $\theta= 1 - 1/ (4C')$ gives the
desired result.

With these two lemmas, we can now complete the proof of the proposition.

Let $\tau_n =0$, and define for $m \le n$, $\tau_m = \inf\{t\ge0\dvtx
N(t) \le3^m\}$. By iterating Lemma~\ref{Linduction} $\lceil\log
_{{\theta}^{-1}} 3\rceil$ times, we see that there exists $C>0$ such that
%
\begin{equation}
\label{decrease} \prob\bigl(\tau_m - \tau_{m-1} \ge C
5^{-m}\bigr) \le C \exp\bigl( - 3^m/C\bigr).
\end{equation}
Let $A_m$ denote the complement of the event above, and let $A = \bigcap_{m=M}^{n} A_m$, where $M$ is a fixed large number. Then on the one
hand, by (\ref{decrease}),
\[
\prob\bigl(A^c\bigr) \le C \exp\bigl( - 3^M/C\bigr)
\]
uniformly in $n$. On the other hand, on the event $A$,
\[
\tau_M \le C5^{-M} + \cdots+ C5^{-n} \le C
5^{-M}. %
\]
Thus, if $\delta>0, \eps>0$ are fixed as in the statement of the
proposition, we choose $M$ large enough that $C 5^{-M} \le\delta$ and
$C \exp( - 3^M/C) \le\eps$. Then picking $k = 3^M$, we obtain
\[
\prob\bigl( N(\delta) >k\bigr) \le\prob( \tau_M \ge\delta) \le\prob
\bigl(A^c\bigr) \le\eps. %
\]
Noting that the choice of $M$ (and thus of $k$) depends only on $\delta
$ and $\eps$ (and not on $n$) completes the proof.
\end{pf}

\subsection{Finitely many coalescing random walks}\label{subsecinvgasket}

Fix $n \ge1$ and let $z_1, \ldots, z_n \in G \times\R$. For $\eta=
2^{-m}$, let $z_1^\eta, \ldots, z_n^\eta$ be space--time points in the
rescaled gasket $\LL_\eta= G_m \times5^{-m} \Z$ such that $z_1^\eta
\to z_1, \ldots, z_n^\eta\to z_n$ as $\eta\to0$. Consider $n$
independent rescaled coalescing simple random walks in $G_m$ started
from $z_1^\eta, \ldots, z_n^\eta$ and making jumps at times in
$5^{-m} \Z$. Let $(Y^\eta_1, \ldots, Y^\eta_n) $, viewed as a
random element of $(\Pi^n, \varrho^n)$ as defined in (\ref
{pathmetric}) with $d=2$, be the collection of $n$ continuous paths
obtained by linearly interpolating the above $n$ coalescing random
walks. Let $(Y_1, \ldots, Y_n)$ be a system of $n$ independent
coalescing Brownian motions started from $z_1, \ldots, z_n$, also
viewed as a random variable in $(\Pi^n, \varrho^n)$.

\begin{proposition} \label{pkpoint}
As $\eta\to0$,
\[
\bigl(Y^\eta_1, \ldots, Y^\eta_n
\bigr) \to(Y_1, \ldots, Y_n) %
\]
in distribution on $(\Pi^n, \varrho^n)$.
\end{proposition}

\begin{pf}
We need the following two facts (\ref{A2pt}) and (\ref{B2pt}) in
the proof.
Let $(W^m_{k5^{-m}})_{k \ge0} $ and $(Z^m_{k5^{-m}})_{k \ge0} $ be
two independent simple random walk on $G_m$ starting at time $0$ and let
\[
\tau^m = \min \bigl\{ k5^{-m} \ge0\dvtx  W^m_{k5^{-m}}
= Z^m_{k5^{-m}} \bigr\}
\]
denote their coalescence time. Then for all $\alpha>0$,
%
\begin{equation}
\label{A2pt} \qquad\limsup_{m \to\infty} \sup_{x,y \in G_m\dvtx  \llvert  x - y\rrvert   \le\eps}
\prob\bigl( \tau^m > \alpha\mid W^m_0 = x,
Z^m_0 = y\bigr) \to0\qquad\mbox{as } \eps \to0.
\end{equation}
Similarly, if $(W_t)_{t \ge0} $ and $(Z_t)_{t \ge0} $ are two
independent Brownian motions on $G$ with coalescence time $\tau$, then
for all $\alpha>0$,
%
\begin{equation}
\label{B2pt} \sup_{x,y \in G\dvtx  \llvert  x - y\rrvert   \le\eps}\prob( \tau> \alpha\mid
W_0 = x, Z_0 = y) \to0\qquad\mbox{as } \eps\to0.
\end{equation}
Let us now prove (\ref{A2pt}). By scaling it is enough to show that
%
\begin{eqnarray}\label{scalingreduction}
\limsup_{m} \sup_{x,y \in G_0\dvtx  \llvert  x - y\rrvert   \le2^m }
\prob\bigl( \tau^0 > 5^{m+M} \mid W^0_0
= x, Z^0_0 = y\bigr) \to0
\nonumber\\[-8pt]\\[-8pt]
\eqntext{\mbox{as } M \to\infty.}
\end{eqnarray}
By Lemma~\ref{Lcoaltimepair}, if $\llvert  x-y\rrvert   \le2^{m}$, then
%
\begin{equation}
\label{hittingprob} \prob\bigl( \tau^0 \le5^{m} \mid
W^0_0 = x, Z^0_0 = y\bigr) \ge c
\end{equation}
for some absolute positive constant $c$. Using the heat kernel upper
bound (\ref{heatub}), it is straightforward to show that given
$\kappa>0$, there exists $K \in\Z_+$ such that for all $m \ge1$,
%
\begin{equation}
\label{eqmaxineqrw} \prob \bigl( \bigl\llvert W^0_{5^m} -
W^0_0\bigr\rrvert > 2^{m+K-2} \bigr) \le\kappa.
\end{equation}
Fix $\eps>0$ and let $L$ be such that $(1-c)^L \le\eps/2 $.
Note that (\ref{hittingprob}) and (\ref{eqmaxineqrw}) with
$\kappa= \eps(4L)^{-1}$ imply that when $m $ large enough, starting
from $x, y \in G_0$ with $\llvert  x - y\rrvert   \le2^{m}$, with probability at least
$c$, the random walks $W^0$ and $Z^0$ either hit each other by time
$5^m$ or else they will be at most $2^m+ 2\cdot2^{m+K-2} \le2^{m+K}$
distance apart at time $5^m$, with probability at least $1 - \eps
(2L)^{-1}$. In the latter case, using the Markov property, again with
probability at least $c$, the random walks $W^0$ and $Z^0$ either hit
each other by the next $5^{m + K}$ \margin{C: I did not quite get why
$K+2$ is needed here, but that's not very important} amount of time, or
else, they will be at most $2^{m+K} + 2 \cdot2^{m+2K-2} \le
2^{m+ 2K}$ distance apart at time $5^m + 5^{m + K}$ with probability at
least $1 - \eps(2L)^{-1}$. Repeat this procedure $L$ times to deduce
that for $M = LK$, the probability in (\ref{scalingreduction}) is
bounded above by $\eps/2 + (1- c)^L \le\eps$ and hence (\ref
{A2pt}) is proved.

To conclude (\ref{B2pt}) from (\ref{A2pt}), use the weak
convergence random walk in the gasket toward Brownian motion, keeping
in mind that
the event $\{ (f, g)\dvtx  f, g \in C(\R_+), f(s) = g(s)$ for some $s \in
[0, \alpha]\}$ is closed in $C(\R_+) \times C(\R_+)$.

Now we are ready to prove Proposition~\ref{pkpoint}. Recall $\eta=
2^{-m}$. Let $X^\eta_1, \ldots, X^\eta_n$ be the continuous paths of
the independent random walks on $G_m$ starting from $z_1^\eta, \ldots,
z_n^\eta$. Similarly, $X_1, \ldots, X_n$ are independent Brownian
motions on $G$ from $z_1, \ldots, z_n$. To keep things simple, we will
assume $n=2$, but the argument can easily be extended to general $n>2$
by induction. By the invariance principle, $X^\eta_i \to X_i$ in
distribution on $(\Pi, \varrho)$. By the Skorokhod representation
theorem, we may assume that $(X^\eta_1, X^\eta_2) \to(X_1, X_2)$
almost surely on $(\Pi^2, \varrho^2)$. Let us show that $\varrho
(Y^\eta_2, Y_2) \to0$ in probability. It is enough to show
for any $R>0$, $\llvert  \hat Y^\eta_2(t) - \hat Y_2(t)\rrvert  _{L^\infty([-R, R])}
\to0$, where for $(\gamma, t_0) \in\Pi$, we denote by $\hat\gamma
$ the continuous function that extends $\gamma$ to all of $\R$ by
setting $\hat\gamma(t) = \gamma(t_0)$ for all $t< t_0$. Define
\[
\eps_\eta= \bigl\llvert \hat X^\eta_1(t) - \hat
X_1(t)\bigr\rrvert _{L^\infty([-R, R])} \vee\bigl\llvert \hat
X^\eta_2(t) - \hat X_2(t)\bigr\rrvert
_{L^\infty([-R, R])}.
\]
We have $\eps_\eta\to0$ almost surely. Let $\tau^\eta= \inf\{ t
\in5^{-m} \Z\dvtx  X^\eta_1(t) = X^\eta_2(t) \}$ and $\tau= \inf\{ t
\in\R\dvtx  X_1(t) = X_2(t) \}$.
Now, we can estimate
\begin{eqnarray*}
\bigl\llvert \hat Y^\eta_2(t) - \hat Y_2(t)
\bigr\rrvert _{L^\infty([-R,
R])} &\le&\eps_\eta + \mathbf{1}_{ \{\tau_\eta> \tau\}}
\sup_{s
\in[ \tau, \tau_\eta] \cap[-R, R] } \bigl\llvert X_2^\eta(s) -
X_1(s) \bigr\rrvert
\\
&&{}+ \mathbf{1}_{ \{\tau_\eta< \tau\}} \sup_{s \in[ \tau_\eta,
\tau] \cap[-R, R] } \bigl\llvert
X_1^\eta(s) - X_2(s) \bigr\rrvert,
\end{eqnarray*}
which is bounded by
%
\begin{eqnarray}\label{eqintermediate}
&& 3\eps_\eta+ \mathbf{1}_{ \{\tau_\eta> \tau\}} \sup
_{s \in[ \tau
, \tau_\eta] \cap[-R, R] } \bigl\llvert X_2(s) - X_1(s)
\bigr\rrvert
\nonumber\\[-8pt]\\[-8pt]\nonumber
&&\qquad{}  + \mathbf{1}_{ \{\tau_\eta< \tau\}} \sup_{s \in[ \tau_\eta,
\tau] \cap[-R, R] } \bigl
\llvert X_1(s) - X_2(s) \bigr\rrvert.
\end{eqnarray}
Note that whenever one of the pairs $(X^\eta_1, X^\eta_2)$ or $(X_1,
X_2)$ coalesces at some time within $[-R, R]$, at that moment the two
processes of the other pair can be at distance at most $2\eps_\eta$
apart from each other. Hence, by (\ref{A2pt}) and (\ref{B2pt}), for
any $\delta>0$,
except for an event with probability at most $ o_\eta(1)$, (\ref
{eqintermediate}) can be bounded by
\[
3\eps_\eta+ \sup_{s \in[ \tau- \delta, \tau+\delta] \cap[-R,
R] } \bigl\llvert
X_1(s) - X_2(s) \bigr\rrvert \le3\eps_\eta+
\osc(X_1- X_2; \delta),
\]
where for a function $\varphi$, its oscillation is defined as $\osc
(\varphi; \delta) = \sup\{\llvert  \varphi(s) - \varphi(t)\rrvert\dvtx  s, t \in
[-R, R], \llvert  s - t\rrvert   \le\delta\}$, and we have used the simple
observation that
$X_1(s) - X_2(s) = (X_1(s) - X_2(s) ) - (X_1(\tau) - X_2(\tau))$.
Since, $X_1$ and $X_2$ are uniformly continuous on $[-R, R]$, almost
surely $\osc(X_1- X_2; \delta) \to0$ as $\delta\to0$.
This proves Proposition~\ref{pkpoint}.
\end{pf}

This completes the proof of Theorem~\ref{tscalinglimitg}.
\end{pf}

\begin{appendix}\label{app}
\section*{Appendix: Connection to the Brownian web of Fontes 
et~al.~\texorpdfstring{\protect\cite{44444}}{[12]}}

\setcounter{theorem}{0}

We now briefly recall from \cite{44444} the construction of the Brownian
web in the space of compact sets of continuous paths. Let $\R^2_c$
denote the completion of the space--time plane $\R^2$ with respect to the
metric
\[
d \bigl( (x_1, t_1),(x_2, t_2)
\bigr) =\bigl\llvert \tanh(t_1) - \tanh(t_2)\bigr\rrvert
\vee\biggl\llvert \frac{\tanh(x_1)}{1 +\llvert   t_1\rrvert  } - \frac{\tanh(x_2)}{1
+\llvert  t_2\rrvert  } \biggr\rrvert.
\]
It is helpful to think of $\R^2_c$ as the continuous image of
$[-\infty, \infty]^2$ under a
map that identifies the lines $[-\infty, \infty] \times\{ \infty\}$
and $[-\infty, \infty] \times\{- \infty\}$ with points $(*, \infty
)$, and $(*, -\infty)$, respectively. Let $\hat\Pi$ be the space of
all continuous paths in $\R^2_c$ with all possible starting times in
$[-\infty, \infty]$. A continuous path $\gamma$ in $\R^2_c$ with
starting point $\sigma_\gamma\in[-\infty, \infty]$, is a mapping
$\gamma\dvtx  [\sigma_\gamma, \infty] \to[-\infty, \infty] \cup\{ *\}
$ such that $\gamma(\infty) = *$, $\gamma(\sigma_\gamma) = *$ if
$\sigma_\gamma= \infty$, and $t \mapsto(\gamma(t), t)$ is
continuous from $[\sigma_\gamma, \infty] \to\R_c^2$. The space
$\hat\Pi$ is equipped with the metric
\begin{eqnarray*}
\hat\varrho \bigl((\gamma_1, \sigma_{\gamma_1}), (
\gamma_2, \sigma _{\gamma_2}) \bigr) &=& \bigl\llvert \tanh(
\sigma_{\gamma_1}) - \tanh (\sigma_{\gamma_2}) \bigr\rrvert
\\
&&{} + \sup_{ t \ge\sigma_{\gamma_1} \wedge\sigma_{\gamma
_2}} \biggl\llvert \frac{\tanh\gamma_1((t \vee\sigma_{\gamma_1}))}{1+
\llvert  t\rrvert  }-
\frac{\tanh(\gamma_2(t \vee\sigma_{\gamma_2}))}{1+ \llvert  t\rrvert  } \biggr\rrvert,
\end{eqnarray*}
which makes $\hat\Pi$ a complete separable metric space. Now define
$\K$ to be the space of compact subsets of $\hat\Pi$ and endow $\K$
with the standard Hausdorff metric $d_\K$. The space $(\K, d_\K)$
again turns out to be a complete separable metric space and the
Brownian web can be defined as a random element in $\K$ via the
following recipe.

Fix a countable ordered set $\D= ( z_1, z_2, \ldots)$, which is dense
in $\R^2$. Let $\cB_n:= \{ W_1^c, \ldots, W_n^c \} $ be the set of
$n$ coalescing Brownian paths starting from the space--time points $z_1,
\ldots, z_n$ as in Section~\ref{SSArratiaflow}, viewed as a random
element in $\K$. In \cite{44444}, the authors showed that almost
surely, $\cB_n$ converges in $\K$ to some random element $\cB_\infty
$, which they called the Brownian web. Its distribution does not depend
on the choice of $\D$. Let $\calW_\infty$ be distributed according
to the coalescing Brownian flow on~$\HH$. It is natural to ask how
these two different objects $\cB_\infty$ and $\calW_\infty$ are related.

For $(\gamma, \sigma_\gamma) \in\hat\Pi$, the notion of a tube
$T$ being crossed (or traversed) by $(\gamma, \sigma_\gamma)$
remains exactly the same. Namely, a tube $T$ is crossed by $(\gamma,
\sigma_\gamma)$ if $ \sigma_\gamma\le t_0$, $(\gamma(t_0), t_0)
\in\partial_0 T$, $(\gamma(t_1), t_1) \in\partial_1 T$ and
$(\gamma(s), s) \in[T]$ for all $s \in(t_0, t_1)$, where $t_0$ and
$t_1$ be the start time and the end time of $T$, respectively. For a
subset\vspace*{1pt} $F$ of $\hat\Pi$, let $\hat\Cr(F)$ denote the set of tubes
in $\TUBE$ which are crossed by at least one path in~$F$. Clearly,
$\hat\Cr(F)$ is always hereditary. The exact same proof of
Lemma~\ref{compactimpliescont} shows that $\hat\Cr(F)$ is closed
in $\TUBE$ if $F$ is compact. This means $\hat\Cr$ maps $\K$ into
$\HH$.

\begin{theorem}\label{thmwebflowconnection}
We have
\[
\hat\Cr( \cB_\infty) \stackrel{d} {=} \calW_\infty.
\]
\end{theorem}

The above theorem says that the law of the coalescing Brownian flow on
$\R$ in the tube topology is nothing but the push-forward of the law
of the Brownian web in the path space. It gives another construction of
the coalescing Brownian flow on~$\R$, though the similar construction
does not work for the gasket due to the absence of appropriate
`Brownian web on the gasket'.
\begin{pf*}{Proof of Theorem~\ref{thmwebflowconnection}}
Let $\calW_n = \calW(z_1, \ldots, z_n) \in\HH$ be as defined in
Section~\ref{SSArratiaflow}. Clearly, $\hat\Cr(\cB_n) = \Cr(\cB
_n) = \calW_n$. We have seen in the proof of Theorem~\ref{tBW} that
\[
\calW_n \stackrel{d} {\to} \calW_\infty\qquad\mbox{in } (
\HH, d_\HH).
\]
Hence, the proof of the theorem is complete by Lemma~\ref{contofPsi}
and the continuous mapping theorem.
\end{pf*}

\begin{lemma} \label{contofPsi}
The map $\hat\Cr\dvtx  (\K, d_\K) \to(\HH, d_\HH)$ is continuous.
\end{lemma}

\begin{pf}
Let $F_n \to F$ in Hausdorff metric $d_\K$. A basis element in the
topology of $\HH$ is of the form $B =\neg\,\boxminus_{T_1} \cap
\cdots\cap\neg\,\boxminus_{T_k} \cap\,\boxminus^{U_1} \cap\cdots
\cap\boxminus^{U_l}$, where $T_1, \ldots, T_k \in\TUBE$ and $U_1,
\ldots, U_l$ are open sets in $\TUBE$. We\vspace*{1pt} need to show that if $H:=
\hat\Cr(F) \in B$, then $H_n:=\hat\Cr(F_n) \in B$ for sufficiently
large $n$.

Note that $H \in\boxminus^{U_i}$ means that there is a $T \in U_i$
such that $T$ is traversed by some $(f, t) \in F$. Since $U_i$ is open,
we can find another tube $T' \in U_i$ such that $T' < T$. Moreover,
there exists $(f_n, t_n) \in F_n$ such that $ \varrho((f_n, t_n), (f,
t) ) \le2 d_\K(F_n, F) \to0$. Clearly, $T'$ is traversed by $(f_n,
t_n)$ for all $n \ge n_0$, and hence, $H_n \in\boxminus^{U_i}$ for
all $n \ge n_0$.

Next, we have to argue that if $T \in\TUBE$ is not traversed by $F$,
then $T \in\TUBE$ is not traversed by $F_n$ for sufficiently large $n$.
Suppose that $T$ is traversed by $(f_n, t_n) \in F_n$ infinitely often,
then it suffices to prove that $T$ is traversed by $F$ as well. For
notational convenience, we will assume that $T$ is traversed by $(f_n,
t_n) \in F_n$ for all~$n$. Using compactness of $F$ and the fact that
$d_\K(F_n, F) \to0$, we can find $(f, t) \in F$ such that $\varrho
((f_n, t_n), (f, t) ) \to0$. It is now easy to check, along the lines
of the proof of Lemma~\ref{compactimpliescont}, that $T$ is
traversed by $(f, t) \in F$. This completes the proof.
\end{pf}

\begin{remark}
By Lemma~\ref{contofPsi}, the pre-image of any event in the tube
topology under the map $\hat\Cr$ is measurable in the path topology.
So, the tube topology is weaker (or coarser) than the path topology.
It would be interesting to find examples of events which are measurable
in the path topology but whose images under $\hat\Cr$ are not
measurable in the tube topology.
\end{remark}
\end{appendix}



%

\printaddresses
\end{document}